\documentclass[12pt,reqno]{amsart}
\usepackage[a4paper,margin=1in]{geometry}

\usepackage{microtype}
\usepackage[cal=euler,scr=boondoxo]{mathalfa}
\usepackage{tikz,tikz-cd}
\usetikzlibrary{decorations}
\usepackage{enumitem}

\usepackage[colorlinks,allcolors=blue]{hyperref}
\usepackage[capitalize]{cleveref}

\input{commands.tex}
\defcite\LQW{li_hilbert}
\defcite\lehn{lehn_chern}
\defcite\GM{goujard_counting}
\defcite\Qin{qin_multiple}

\crefname{equation}{}{}
\crefname{conjecture}{Conjecture}{Conjectures}

\opr\iend{\mathbf{End}}
\opr\hom{\mathbf{Hom}}
\def\hb#1#2{#1^{[#2]}}
\def\pr#1#2{#1^{(#2)}}

\def\P{\opn{Pair}}
\def\Pair(#1){\P_{#1}}
\def\d{\mathbf d}
\def\grd{\mathbf d} 

\operl{BD,Amp,CT,QS}
\oper{str} 
\def\hh{\what H}

\def\pN{\bN_*} 
\def\Zs{\bZ_*}
\def\Ns{\pN}

\def\fat{\fa^\ts}
\opr\fch{\mathfrak{ch}} 

\def\N{\sN} 
\def\Ho{\cH} 

\def\Z{\bZ}

\oper{Contr}
\oper{Fib}

\newlist{clm-l}{enumerate}{1} 
\setlist[clm-l]{
  label=\bfseries(\arabic*),
  itemindent=*,
  leftmargin=0em,
  labelsep=0.5em,
}

\def\bb#1#2{\pser{#1;#2}} 

\def\Pg{\cP^{\mathrm g}}

\def\ld{\mathfrak d} 
\def\W{\cW}
\def\bone{\mathbf1}
\def\mset#1{\set{\hspace{-3pt}\set{#1}\hspace{-3pt}}}
\def\mset{\sbr}

\def\ct{\ang}
\def\hB{\what B}
\def\hD{\what D}

\def\bde{\boldsymbol\vartheta} 

\def\hB{\bfB}
\def\hD{\bfD}
\def\hg{\hat g}
\def\tG{\wtl\fG}
\oper{ev}

\def\M{\cM} 

\def\PM{\mathcal{PM}} 
\def\P{\PM} 
\def\pip{\pi^{(2)}} 
\def\pis{\pi^{(1)}} 
\def\ph#1{(#1;q)_\infty}

\def\prods{\prod\nolimits}
\def\sums{\sum\nolimits}

\def\pp#1{[#1]_+} 
\def\pn#1{[#1]_-} 
\oper{QM}

\def\Zq{\cZ_q}
\def\RR{\mathsf R} 
\def\BB{\mathsf B}
\def\zz{w}
\def\hPhi{\what\Phi}
\def\odd{\mathrm{odd}}
\def\bm{\mathbf m}
\oper{Or}

\def\Sa{S}
\def\ra{\rho} 
\def\rb{\rho_2} 

\begin{document}
\title{Integrals over Hilbert schemes and multiple $q$-zeta values}

\author{Sergey Mozgovoy}

\address{School of Mathematics, Trinity College Dublin, Dublin 2, Ireland
\newline\indent
Hamilton Mathematics Institute, Dublin 2, Ireland}

\email{mozgovoy@maths.tcd.ie}

\begin{abstract}
We establish a version of Okounkov's conjecture on multiple $q$-zeta values arising from integrals over Hilbert schemes of points on smooth projective surfaces. We prove that normalized generating series of integrals involving products of Chern characters of tautological bundles belong to the algebra $\mathcal Z_q$ of multiple $q$-zeta values, without restrictions on the canonical class of the surface. The proof combines vertex-operator methods with a reduction to graph sums, whose arithmetic structure is governed by the cohomology of the surface. We also derive a closed universal formula for the Chern character operators acting on the cohomology of Hilbert schemes by solving a system of differential equations arising from Lehn's commutator relations.
\end{abstract}

\maketitle

\section*{Introduction}

Generating series of intersection numbers on moduli spaces often exhibit arithmetic structures that are not visible in their individual coefficients.
For Hilbert schemes of points on surfaces, this phenomenon takes a particularly striking form: the geometry of tautological bundles is expected to produce $q$-series belonging to distinguished algebras of multiple $q$-zeta values.
This conjecture belongs to a broader circle of ideas relating enumerative geometry, vertex operators, and the modular properties of partition functions.
The purpose of this paper is to establish a multiple $q$-zeta value property for a large class of Hilbert-scheme integrals, without imposing any condition on the canonical class of the surface.

Let $X$ be a smooth projective complex surface, $\hb Xn$ its Hilbert scheme of $n$ points, and $\hb En$ the tautological vector bundle associated with a vector bundle $E$ on $X$.
For vector bundles $E_1,\ldots,E_r$ and integers
$k_1,\ldots,k_r\ge0$, consider the $q$-bracket
\begin{equation}\label{intro:bracket}
\ang{\ch_{k_1}^{E_1}\dots \ch_{k_r}^{E_r}}
=\sum_{n\ge0}q^n\int_{\hb Xn}
\ch_{k_1}(\hb{E_1}n)\dots \ch_{k_r}(\hb{E_r}n)
\cdot \fe_t(T_{\hb Xn}),
\end{equation}
where $\fe_t(E)=\fe(E\ts e^t)=\prod_i(t+x_i)$ for a vector bundle $E$
with the Chern roots $x_i$.
Define the reduced (or normalized) $q$-bracket $\ang{f}'=\ang{f}/\ang{}$,
where the empty bracket is the generating function of Euler characteristics,
\[
\ang{}
=\sum_{n\ge0}\chi(\hb Xn)q^n
=\prod_{i\ge1}(1-q^i)^{-\chi(X)}.
\]

Okounkov conjectured \cite{okounkov_hilbert} that,
for tautological line bundles, the reduced $q$-brackets belong
to the algebra $\cZ_{q,1}$ of multiple $q$-zeta values,
spanned by the elements
\eq{\label{Z:int}
Z(f_1\ts\dots\ts f_r)
=\sum_{n_1>\dots>n_r>0}\prod_{i=1}^{r}f_i(q^{n_i})\in\bQ\pser q,
}
where $f_i\in \RR_1=\sets{g/(1-q)^k}{g\in\bQ[q],\, \deg g<k}$
satisfy $f_i(0)=0$,
see \eg \cite{bachmann_dimension,hirose_unified,mozgovoy_multiple}.
This prediction gives a remarkably restrictive arithmetic description of intersection numbers on the Hilbert schemes~$\hb Xn$.

A related conjecture of Qin \cite{qin_quick,alhwaimel_hilbert}
asserts quasi-modularity when the canonical class of $X$ is numerically trivial.
This was recently proved in \cite{alekseev_qins,hongminh_quasia}
(see also \cite{carlsson_vertex,carlsson_vertexa,alhwaimel_qins,alhwaimel_equivariant}).

Our main result concerns the algebra $\Zq$ of multiple $q$-zeta values, spanned by the elements~\cref{Z:int},
where $f_i$ belong to
\[
\RR=\sets{g/(1-q)^k}{g\in\bQ[q],\, \deg g\le k}
=\bQ\sbr{\tfrac q{1-q}}\]
and satisfy $f_i(0)=0$, see \eg \cite{bachmann_dimension,hirose_unified,mozgovoy_multiple}.
This algebra coincides with the algebra of \idef{bi-brackets}
\cite{bachmann_algebraa,bachmann_dimension}
and contains the algebra $\QM$ of quasi-modular forms.
We prove the following result.

\begin{theorem}\label{intro:main}
Let $X$ be a smooth projective surface,
$E_1,\ldots,E_r$ be vector bundles on $X$ and
$k_1,\ldots,k_r\ge0$.
Then the reduced $q$-bracket satisfies
\[
\ang{\ch_{k_1}^{E_1}\dots\ch_{k_r}^{E_r}}'\in\Zq[t].
\]
\end{theorem}

For simplicity, we assume that $H^\odd(X,\bQ)=0$ and work with ordinary vector spaces rather than superspaces;
the arguments extend to arbitrary surfaces using the superalgebra formalism (\cf \cite{hongminh_quasia}).
The cohomology $H=H^*(X,\bQ)$ is then a Frobenius algebra, with associated Heisenberg algebra $\what H$.
By the results of Nakajima and Grojnowski \cite{nakajima_heisenberg,grojnowski_instantons}, the space
$\bH=\bop_{n\ge0}H^*(X^{[n]},\bQ)$
is canonically realized as a Fock representation of~$\what H$.
The Carlsson--Okounkov vertex operator
\cite{carlsson_vertex,carlsson_exts} converts the reduced brackets into normalized $q$-traces of products of Chern character operators
on $\bH$ \cite{lehn_chern}.

The geometry is thereby encoded in the Chern character operators
and closely related tautological operators
\cite{li_vertex,li_hilbert}.
These operators admit normally ordered Heisenberg expansions
involving universal contributions of $1$, the Euler class $e_X$, the canonical class $K_X$, and $K_X^2$ \cite{li_hilbert}.
We determine their coefficients by translating Lehn's commutator
relations~\cite{lehn_chern} into a system of differential equations
and solving this system.
The resulting closed formulas are finite sums indexed by cuts of the Heisenberg modes.
A different formula for these coefficients was recently obtained in \cite{li_canonical} by solving a system of recurrence relations.

The main challenge is to understand the arithmetic structure of the resulting generating series. A twisted Wick formula translates the geometric problem into graph sums, whose complexity is strongly constrained by the cohomology of a surface. Contributions involving more than one independent cycle vanish, while those with a single cycle are quasi-modular. The remaining tree contributions encode the dependence on the canonical class and require a more general arithmetic framework. Crucially, the dimension of the surface restricts their complexity, allowing us to establish their membership in the algebra of bi-brackets through Fourier analysis and Lambert series. An essential ingredient is a finiteness property inherited from the original geometric expressions, which ensures the cancellation of potentially divergent terms. In this way, our approach extends the constant-term methods underlying the quasi-modularity of Feynman-like graph sums \cite{goujard_counting} to the more general setting of Hilbert-scheme intersection theory.

The paper is organized as follows.
In \cref{fock:sec}, we recall the theory of Heisenberg algebras and Fock spaces associated with Frobenius algebras, including graph contractions.
We also introduce a polynomial realization of normally ordered Heisenberg monomials.
In \cref{hilb:sec}, we review Nakajima operators, Chern character operators, tautological operators, and Lehn's boundary operator, together with their fundamental commutator relations.
In \cref{diff:sec}, we translate these relations into a system of differential equations and solve it to obtain closed universal formulas for the Chern character operators.
In \cref{QM:sec}, we develop the theory of quasi-elliptic functions and their Fourier coefficients, establish quasi-modularity results for graph sums, and describe the structure of Fourier coefficients in terms of a hyperplane arrangement.
In \cref{mqzv:sec}, we introduce multiple $q$-zeta values and bi-brackets and prove the summation results needed for our applications, including a general theorem for nested families.
In \cref{wick:sec}, we establish Wick and twisted Wick formulas for normalized $q$-traces of Heisenberg operators.
Finally, in \cref{ok:sec}, we use the Carlsson--Okounkov vertex operator and the twisted Wick formula to reduce the main theorem to graph sums. We analyze the contributions of graphs with cycles and reduce the remaining tree contributions to a framed-tree summation theorem, which we prove using the Fourier coefficient and bi-bracket summation results developed earlier.

\subsection*{AI Declaration}
An LLM was used for exploratory discussions and editorial assistance
during the preparation of this paper.

\subsection*{Conventions}
We write $\bN=\set{0,1,2,\dots}$,
$\Ns=\bN\ms\set0$, $\Zs=\bZ\ms\set{0}$,
and $[n]=\set{1,\dots,n}$ for $n\ge0$.
For a series $f(z)=\sum_n f_nz^n$, we write $[z^n]f=f_n$ for the coefficient of $z^n$ in $f$.

\section{Fock spaces}
\label{fock:sec}
In this section we will discuss the general theory of Heisenberg algebras and Fock spaces associated with quadratic spaces
and with Frobenius algebras.
This will be used later in the case of the Frobenius algebra
$H=H^*(X,\bQ)$ associated with a smooth projective surface over $\bC$.
We separate the results about Fock spaces that depend only on a Frobenius algebra, from the results that depend on a surface.
However, some structures and results that depend only on a Frobenius algebra require a geometric motivation, and therefore are postponed to later sections.

\subsection{Heisenberg algebra and Fock space}
A quadratic space over a field $\bk$ of characteristic $\ne2$
is a \fdd vector space $H$ equipped with a non-degenerate symmetric bilinear form $(-,-)$.
The Heisenberg algebra of $H$ is the graded Lie algebra
\[\hh=H[t^{\pm1}]\oplus\bk\bc\]
with the central element $\bc$ and the Lie bracket
\[
[a_m,b_n]=m\de_{m,-n}(a,b)\bc,\qquad a_m=at^m,\,b_n=bt^n,
\]
for $a,b\in H$ and $m,n\in\bZ$.
The degree on $H$ is defined by $\deg a_n=-n$ and $\deg\bc=0$.
Let $\hh_+=H[t]$ and $\hh_-=t\i H[t\i]$.
The Fock space of $H$ is the $\hh$-module
\[
V=V(H)=\Ind^{\hh}_{\hh_+\oplus\bk\bc}\bk\iso S(\hh_-),
\]
where $\hh_+$ acts as zero and $\bc$ acts as identity on $\bk$.
It inherits the grading $V=\bop_{n\ge0}V_n$ such that $\dim V_0=1$ and
$\dim V_n<\infty$ for all $n$.
The element $\one=\ket0\in V_0$ corresponding to $1\in\bk$ is called the vacuum.

By abuse of notation,
we denote by $a_n=at^n$ the element of~$\hh$ (a Heisenberg mode)
and its image in $\End(V)$ (a Heisenberg operator).
The operators $a_n$ for $n\ge0$ are called annihilators and
the operators $a_n$ for $n<0$ are called creators.
For $n\in\bZ$, consider the linear map
\[\fa_n:H\to\End(V),\qquad
a\mto a_n.
\]

Given two $\bZ$-graded vector spaces $V,W$, define the internal Hom
\[\hom(V,W)=\bop_{d\in\bZ}\hom_d(V,W),\qquad
\hom_d(V,W)=\prod_{i\in\bZ}\Hom(V_i,W_{i+d}).
\]
Then $\iend(V)=\hom(V,V)$ is a graded algebra.
In particular, for the Fock space $V=V(H)$,
the representation map $\hh\to\iend(V)$ preserves degrees
and we have $a_n\in\iend_{-n}(V)$.

\subsection{Frobenius algebras}
Later we will study Fock spaces associated with commutative Frobenius algebras.
A Frobenius algebra over a field $\bk$ is a
finite-dimensional (associative, unital) algebra
$(H,\mu,\eta)$ equipped with a linear map $\eps:H\to \bk$,
called the trace or counit, such that the pairing
$(a,b)=\eps(ab)$ on $H$ is non-degenerate
(see \eg \cite{kock_frobenius}).
It has a dual coalgebra structure $(H,\de,\eps)$, where
$\de$ is defined by
\eq{\label{frob:pairing}
(\de(a),b\ts c)=(a,cb),\qquad
(a\ts b,x\ts y)=(a,x)(b,y).
}
The algebra and coalgebra structures satisfy
\begin{equation}\label{frob:eq}
\de\mu=(\mu\ts\id)(\id\ts\de)
=(\id\ts\mu)(\de\ts\id).
\end{equation}
Conversely, given algebra and coalgebra structures on $H$
satisfying \cref{frob:eq}, the pairing
$(a,b)=\eps(ab)$ is automatically non-degenerate
and \cref{frob:pairing} holds.

The element $e=\mu\de(1)\in H$ is called the Euler element.
For every $a\in H$, we have
\eq{\mu\de(a)=\mu\de(1\cdot a)
=\mu(\id\ts\mu)(\de 1\ts a)
=\mu(\mu\ts\id)(\de 1\ts a)=ea.
}

Define $\de_r:H\to H^{\ts r}$ by $\de_0=\eps$, $\de_1=\id$,
and generally $\de_{r+1}=(\de_r\ts\id)\de$ for $r\ge0$.
Similarly, define $\mu_r:H^{\ts r}\to H$.
For $1\le i<j\le r$, define the contraction
\[\fc_{ij}:H^{\ts r}\to H^{\ts(r-2)},\qquad a_1\ts\dots\ts a_r
\mto(a_i,a_j)a_1\ts\dots\ts\what a_i\ts\dots\ts
\what a_j\ts\dots\ts a_r.
\]

In what follows, we assume that $H$ is a commutative (and cocommutative) Frobenius algebra.

\begin{lemma}[Frobenius contraction]
\label{frob contr}
Let $a,b\in H$.
For $r,s\ge1$, we have
\eq{\label{frob1}
\fc_{ij}\rbr{\de_r(a)\ts\de_s(b)}=\de_{r+s-2}(ab),
\qquad 1\le i\le r<j\le r+s.
}
For $r\ge2$, we have
\eq{\label{frob2}
\fc_{ij}(\de_r(a))=\de_{r-2}(ea),\qquad 1\le i<j\le r.
}
\end{lemma}
\begin{proof}
\clm1
Assume for simplicity that $\de(a)=a_1\ts a_2$ and $\de(b)=b_1\ts b_2$.
Let
$c=c'\ts c''\in H^{\ts r-1}\ts H^{\ts s-1}$
and let $x=\mu_{r-1}(c')$, $y=\mu_{s-1}(c'')$.
Then
\begin{multline*}
\rbig{\fc_{ij}\rbr{\de_r(a)\ts\de_s(b)},c}
=(a_1,b_1)\rbig{\de_{r-1}(a_2)\ts\de_{s-1}(b_2),c'\ts c''}
=(a_1,b_1)(a_2,x)(b_2,y)\\
=(a,b_1x)(b_2,y)
=(b_1,ax)(b_2,y)=(b,axy)=(\de_{r+s-2}(ab),c).
\end{multline*}
Therefore $\fc_{ij}\rbr{\de_r(a)\ts\de_s(b)}=\de_{r+s-2}(ab)$.

\clm2
Applying the previous result to $\de_r(a)=(\de_{r-1}\ts \de_1)\de(a)$,
we obtain $\fc_{ij}(\de_r(a))=\de_{r-2}(\mu\de(a))=\de_{r-2}(ea)$.
\end{proof}

%
%

\begin{lemma}[Frobenius graph contraction]
\label{lem:frobenius-graph}
Let $G=(V,E)$ be a connected graph and let $a=(a_v)_{v\in V}\in H^V$.
For each vertex $v\in V$,
let $n_v$ be its valence and let the factors of $\de_{n_v}(a_v)$
be placed at the incident half-edges (flags).
Let $\cF_G(a)$ be obtained by applying Frobenius pairing to every edge and taking the product.
Then
\begin{equation}\label{eq:frobenius-graph}
\cF_G(a)=\eps\rbr{e^{b_1(G)}\prod*_{v\in V}a_v}.
\end{equation}
\end{lemma}

\begin{proof}
Choose a spanning tree $T$ of $G$.
The Frobenius identity \cref{frob1}
shows that contracting an edge joining two distinct
vertices merges them into one vertex, whose label
is the product of their labels.
Contracting all edges of $T$ therefore produces
one vertex labeled by $c=\prod_{v\in V}a_v$.
The remaining edges become loops.
By the Frobenius identity \cref{frob2}
each loop contraction multiplies the vertex label
by the Euler element $e$.
After contracting all loops, we obtain
$\de_0(e^{b_1(G)}c)=\eps(e^{b_1(G)}c)$.
\end{proof}

\begin{example}
Let $G$ be a graph with one vertex and one loop
and let $\de(a)=a_1\ts a_2$.
The above result implies that $\eps(a_1a_2)=\eps(ea)$.
This follows directly from the fact that $a_1a_2=\mu\de(a)=ea$.

Let $G$ be a graph with two vertices and two edges between them
and let $\de(a)=a_1\ts a_2$ and $\de(b)=b_1\ts b_2$.
The above result implies that
\[
\eps(a_1b_1)\eps(a_2b_2)=\eps(eab).
\]
This follows directly from the fact that
$(\de(a),\de(b))=(\mu\de(a),b)=(ea,b)=\eps(eab)$.
\end{example}

\subsection{Heisenberg monomials}
Let $H$ be a commutative Frobenius algebra.
The pairing $(a,b)=\eps(ab)$ on $H$ is non-degenerate,
hence we may consider the corresponding Heisenberg Lie algebra $\hh$ and the Fock space $V=V(H)$. 
Consider the graded algebra $\iend(V)$ 
and let
\[\N:\iend(V)^{\ts r}\to \iend(V)\]
be the normally ordered product
(see \eg \cite{hongminh_quasia};
the factors are arranged in decreasing order of degree).
We also write $\nop{f_1\dots f_r}=\N(f_1\ts\dots\ts f_r)$.
For $\bn\in\bZ^r$, define
\[\fat_\bn=\fa_{n_1}\ts\dots\ts\fa_{n_r}:
H^{\ts r}\to\iend(V)^{\ts r}\]
and define the (normally ordered) Heisenberg monomial
\[\fa_\bn=\sN\fat_\bn\de_r:H\to\iend(V).
\]
Since $\de$ is cocommutative,
the map $\fa_\bn$ depends just on the orbit
$\la\in\bZ^r/S_r$ of $\bn$, and we denote it by $\fa_\la$.
We also denote~$\fa_\bn$ by $\nop{\fa_{n_1}\dots\fa_{n_r}}\de_r$.

\begin{remark}[Multisets and partitions]
A finite multiset on a set $X$
is an element of $\sS(X)=\bigsqcup_{r\ge0}X^r/S_r$.
For example, for $\Ns=\bN\ms\set0$,
the elements of $\sS(\Ns)$ are partitions.
There is a canonical identification
\[\sS(X)\iso\bN^{(X)}=\sets{m:X\to\bN}{\#\supp(m)<\infty},\]
where a multiset $\la=\mset{\la_1,\dots,\la_r}\in X^r/S_r$
corresponds to $m=m(\la)$ with
\[m_x=\#\sets{i\in[r]}{\la_i=x},\qquad x\in X.\]
We write $\la=(x^{m_x})_{x\in X}$ and define
the length $\ell(\la)=\n m=\sum*_{x\in X} m_x$
and $\la!=\prod_{x\in X}m_x!$.
If $f:X\to A$ is a map to a commutative monoid,
we define
\eq{\sum_{x\in\la}f(x)
:=\sum_{x\in X}m_x f(x).
}
If $X$ is a ring, we define
\eq{
\n\la:=\sum_{x\in X}m_x x,\qquad
p_k(\la):=\sum_{x\in\la}x^k=\sum_{x\in X}m_x x^k,\qquad k\ge1.
}
\end{remark}

For $\Zs=\bZ\ms\set0$,
the elements of $\Pg=\sS(\Zs)$ will be called generalized partitions.
Let
\eq{\label{Pg}
\Pg_{r,n}=\sets{\la\in\Pg}{\ell(\la)=r,\,\n\la=n},\qquad
r\ge0,\, n\in\bZ.
}
We call $\n\la$ the charge (or weight) of $\la$.

\begin{lemma}
\label{a-comm1}
Let $\la\in\Pg$ and let $k\ne0$.
Then
\begin{equation}
[\fa_\la(a),\fa_{-k}(b)]
=k\,m_k(\la)\,\fa_{\la-\set k}(ab),
\label{LQW:basic-comm}
\end{equation}
where the right-hand side is zero if $m_k(\la)=0$.
\end{lemma}
\begin{proof}
Only a factor $\fa_k$ can have a nonzero commutator with
$\fa_{-k}(b)$.
The Heisenberg relation gives
$[\fa_k(c),\fa_{-k}(b)]=k(c,b)\id$ and there are $m_k(\la)$ such factors.
Contracting the corresponding tensor factor of $\de_{\ell(\la)}(a)$
against $b$ gives $\de_{\ell(\la)-1}(ab)$.
This implies \cref{LQW:basic-comm}.
\end{proof}

For $r\ge0$ and $n\in\bZ$, define
\eq{\label{Wrn}
\W^{r}_n(a)
:=\sum_{\ov{\ell(\la)=r}{\n\la=n}}\frac{1}{\la!}\fa_\la(a)
=\frac1{r!}\sum_{\ov{\bn\in\Zs^{r}}{\n\bn=n}}\fa_{\bn}(a).
}
In particular, $\W^0_n(a)=\de_{n,0}\eps(a)\id_V$
and $\W^1_n(a)=\fa_n(a)$.
Define
\eq{\label{oper L}
\fL_n(a):=\W^2_n(a)
=\frac12\sum_{i+j=n}\fa_{i,j}(a)
=\frac12\sum_{i+j=n}\nop{\fa_i\fa_j}\de(a).
}
Using $\fa(a,z)=\fa(z)(a):=\sum_{n\in\bZ}\fa_n(a)z^{-n-1}$,
we obtain
\[\W^{r}(a,z):=\sum_{n\in\bZ}\W^{r}_n(a)z^{-n-r}
=\frac1{r!}\nop{\fa(z)^{r}}\de_{r}(a).
\]

\begin{corollary}\label{a-comm2}
For $r\ge1$ and $m,n\in\bZ$, we have
\[[\W^r_m(a),\fa_n(b)]=-n \W^{r-1}_{m+n}(ab).\]
\end{corollary}
\begin{proof}
For $\mu=\la-\set{-n}$,
we have $m_{-n}(\la)/\la!=1/\mu!$.
\end{proof}

\begin{corollary}
Let $T:=\fL_{-1}(1)=\oh\sum_{i+j=-1}\fa_{i,j}(1)$.
Then $[T,\fa_n(a)]=-n\fa_{n-1}(a)$ and
\[\fa_{-n-1}(a)=\frac{1}{n!}\ad^n_{T}(\fa_{-1}(a)),\qquad
n\ge0.\]
\end{corollary}
\begin{proof}
We have
$[T,\fa_{n}(a)]
=[\W^2_{-1}(1),\fa_{n}(a)]
=-n\fa_{n-1}(a).$
\end{proof}

\begin{remark}
The operator $\d=\fL_0(1)$ is the degree operator,
since $\d\one=0$ and $[\d,\fa_n(a)]=[\W^2_0(1),\fa_n(a)]
=-n\fa_n(a)$.
\end{remark}

\subsection{Polynomial realization}
In this section, we introduce a commutative algebra $\cA$ that encodes normally ordered Heisenberg monomials.
We show that commutators with certain natural operators on the Fock space correspond to differential operators on $\cA$,
allowing us to characterize certain operators on the Fock space as solutions of differential equations.

Consider the algebra
\eq{\cA
=\bop_{r\ge0}\bop_{n\in\bZ}\cA_{r,n}
\sbs\bQ\pser{x_i\col i\in\Zs},\qquad
\cA_{r,n}=\prod_{\nu\in\Pg_{r,n}}\bQ x_\nu,}
where $x_\nu=x^{m(\nu)}=\prod_{i\ne0}x_i^{m_i(\nu)}$.
Let $\dd_i=\dd/\dd x_i$ for $i\in\Zs$.
We also define $x_0=0$ and $\dd_0=0$.

For $a\in H$, define the linear map
\eq{\label{N-A}
\N_a:\cA_H=\cA\ts H\to \End(V),\qquad
x_\nu\ts b\mto \fa_\nu (ba).
}

\begin{remark}
We usually denote $f\ts b\in \cA_H$ by $bf$.
For a linear operator $L$ on $\cA$,
we define the linear operator $bL:=m_bL=Lm_b$ on $\cA_H$,
where $m_b(a)=ba$ for $a\in H$.
\end{remark}

\begin{remark}
In what follows, we identify
$\cA_{r,n}$ with $\cF_{r,n}=\Map(\Pg_{r,n},\bQ)$
via the isomorphism
\[\cF_{r,n}\isoto\cA_{r,n},\qquad
f\mto \bar f=\sum_{\nu\in\Pg_{r,n}} \frac{f(\nu)}{\nu!}x_\nu.
\]
Thus, all operators on $\cA_{r,n}$ considered below
induce operators on $\cF_{r,n}$.
For $a\in H$, define
\eq{
\N'_a:\cF_{r,n}\to\End(V),\qquad
f\mto \sum_{\nu\in\Pg_{r,n}}\frac{f(\nu)}{\nu!}\fa_\nu(a).
}
Then $\N'_a(f)=\N_a(\bar f)$ for $f\in\cF_{r,n}$.
Note that $\N'_a(\bone_{r,n})=\W^{r}_n(a)$ \cref{Wrn}
for the identity function~$\bone_{r,n}$.
\end{remark}

The motivation for the linear map \cref{N-A} is that certain adjoint operators on $\End(V)$ correspond to differential operators on $\cA$.
As a result, certain identities between operators on $\End(V)$
correspond to systems of differential equations on $\cA$.
The next result illustrates this philosophy.

\begin{lemma}\label{N-comm}
For $a,b\in H$ and $f\in\cA$, we have
\begin{gather}
[\N_a(f),\fa_{-1}(b)]=\N_{ab}(\dd_1 f),\\
[T,\N_a(f)]=\N_a(\cT f),\qquad
\cT=-\sum_{n\in\bZ}nx_{n-1}\dd_n. \label{trans}
\end{gather}
\end{lemma}
\begin{proof}
\clm1
For $f=x_\nu\ts c$, we have
(see \cref{a-comm1})
\[
[\N_a(f),\fa_{-1}(b)]=[\fa_\nu(c a),\fa_{-1}(b)]
=m_1(\nu)\fa_{\nu\ms\set{1}}(c ab)=\N_{ab}(\dd_1 f).
\]
\clm2
We have $[T,\fa_n(a)]=-n\fa_{n-1}(a)$.
Therefore, for $f=x_\nu\ts c$, we have
\[
[T,\N_a(f)]=[T,\fa_\nu(c a)]
=-\sum*_n n m_n(\nu)\fa_{\nu\ms\set{n}\cup\set{n-1}}(c a)
=\N_a(\cT x_\nu\ts c).
\qedhere
\]
\end{proof}

\begin{remark}
We have
\[[\cT,\dd_n]=\sum_{k\in\bZ}k\dd_n(x_{k-1})\dd_k=(n+1)\dd_{n+1},\qquad
n\ne0.
\]
\end{remark}

\section{Hilbert schemes and tautological operators}
\label{hilb:sec}
We recall the realization of the cohomology of Hilbert schemes as a Fock space and the operators arising from tautological characteristic classes.
This allows us to translate geometric questions about intersection numbers into algebraic questions about Heisenberg operators and their commutation relations.
For more information on Hilbert schemes see \eg
\cite{nakajima_lectures,lehn_lectures,qin_hilbert}.

\subsection{Nakajima operators}
Let $X$ be a smooth projective surface over $\bC$.
For simplicity, we assume that $H^\odd(X,\bQ)=0$.
Then $H=H^*(X,\bQ)$ is equipped with the Frobenius algebra structure, with the usual cup product and the pairing
\[
(a,b)=-\int_X ab.
\]
The coproduct and the counit are given by
\[\de=-\De_*:H\to H\ts H,\qquad \eps=-\int_X:H\to\bQ,\]
where $\De:X\to X\xx X$ is the diagonal.

\begin{remark}
The Euler element of $H$ is
$e=\mu\de(1_X)=-\De^*\De_*(1_X)=-c_2(T_X)=-e_X$.
\end{remark}

The corresponding Fock space $V=V(H)$ is canonically isomorphic to \[\bH=\bop_{n\ge0}\bH_n,\qquad \bH_n=H^*(\hb Xn,\bQ),\]
as a module over the Heisenberg algebra $\hh$.
More precisely, let $\one\in\bH_0=H^*(\hb X0,\bQ)\iso\bQ$,
called the \idef{vacuum}.
For $n>0$ and $a\in H$, define the operator
$\fa_{-n}(a):\bH\to\bH$ of degree $n$,
called the \idef{Nakajima operator}, as follows.
For $m\ge0$, let
\[\hb Q{m+n,m}
=\sets{(\ksi,x,\eta)\in\hb X{m+n}\xx X\xx \hb Xm}
{\ksi\sps\eta,\,\supp(\ksi-\eta)=\set x}\]
and let $p_1,p_2,p_3$ denote the corresponding projections.
Then
\[\fa_{-n}(a):\bH_m\to\bH_{m+n},\qquad
u\mto p_{1*}\rbr{[\hb Q{m+n,m}]\cdot p_2^*a\cdot p_3^*u}.
\]
Let $\fa_{n}(a)=(-1)^n\fa_{-n}(a)^*$,
where $f^*:\bH_{m+n}\to \bH_m$ for $f:\bH_m\to\bH_{m+n}$
is the adjoint \wrt the intersection pairings on $\bH_k$
\[(fu,v)=(-1)^{p(f)p(u)}(u,f^*v),\qquad u\in \bH_m,\, v\in\bH_{m+n}.\]
Let $\fa_0(a)=0$.
This defines linear maps
\[\fa_n:H\to\End(\bH),\qquad n\in\bZ.\]

The following result is proved in
\cite{nakajima_heisenberg,grojnowski_instantons}
(see also \cite{ellingsrud_intersection,lehn_chern}).

\begin{theorem}
The Nakajima operators satisfy the commutation relations
\[[\fa_m(a),\fa_n(b)]=m\de_{m,-n}(a,b)\cdot \id_\bH,\]
hence define an action of the Heisenberg algebra $\hh$ on $\bH$.
Moreover, for the Fock space~$V$,
there is a unique $\hh$-module isomorphism
\[V\isoto \bH\]
sending the vacuum vector of $V$ to the vacuum vector of $\bH$.
In this way we identify the Nakajima operators $\fa_n(a)$
and the Heisenberg operators $a_n$ on the Fock space.
\end{theorem}

\subsection{Tautological bundles and tautological operators}
\label{sec:taut bun}
Given a vector bundle $E$ on~$X$,
its tautological bundle on $\hb Xn$ is
\eq{
\hb En=p_{1*}(\cO_{Z_n}\ts p_2^*E),
}
where $Z_n\sbs\hb Xn\xx X$ is the universal subscheme and
$p_1,p_2$ are projections.
This is a vector bundle of rank $n \rk(E)$.
Consider the \idef{Chern character operator} \lehn
\eq{
\fch_k(E):\bH\to\bH,\qquad \bH_n\ni u\mto \ch_k(\hb En)\cdot u\in\bH_n.
}
This defines a linear map $\fch_k:K(X)\to\iend_0(\bH)$.

Similarly, for $a\in H$, define the \idef{tautological classes}
\cite{li_vertex}
\eq{G(a,n)=\hb an
=p_{1*}\rbig{\ch(\cO_{Z_n})\cdot p_2^*\rbig{a\td(X)}}\in\bH_n.
}
By Grothendieck–Riemann–Roch,
\eq{\ch(\hb En)=G(\ch(E),n).}
For $a\in H^i$, let $G_k(a,n)$ denote the component of $G(a,n)$ in $\bH_n^{2k+i}$.
Then
\eq{\label{ch-G}
\ch_k(\hb En)=\sum_{i=0}^k G_{k-i}(\ch_i(E),n).
}
The classes $G_k(a,n)$ induce the \idef{tautological operator}
\cite{li_vertex}
\eq{\fG_k(a):\bH\to\bH,\qquad \bH_n\ni u\mto G_k(a,n)\cdot u\in\bH_n.
}
This defines a linear map $\fG_k:H\to\iend_0(\bH)$
such that
$\fch_k=\sum_{i=0}^k\fG_{k-i}\circ\ch_i$.
For example, $\fG_0(1_X)=\grd$ is the degree operator,
since $G_0(1_X,n)=\ch_0(\hb{\cO_X}n)=n$.

An algorithm to determine the Chern character operators $\fch_k(E)\in\End(\bH_n)$ for small $n$ was proposed in \cite[\S4.1]{lehn_chern}.
This idea was used in \cite[Theorem 4.7]{li_hilbert}
to prove that there exist universal functions $g_\ep(\la)$ on generalized partitions,
for $\ep\in\cE=\set{1,e,K,K^2}$ (where $K=K_X\in H^2$ is the canonical class and $e=\mu\de(1)=-e_X$), such that
\eq{ \label{LQW:general}
\fG_k(a)=(-1)^k
\sum_{\ep\in\cE}\sum_{\ov{\n\la=0}{\ell(\la)=k+2-\n\ep/2}}
\frac{g_\ep(\la)}{\la!}\fa_\la(\ep a).
}
Moreover, $g_1$ and $g_e$ were explicitly determined in \LQW
and it was proved that
\eq{\fG_0(a)=\fL_0(a)=\sum_{n>0}\fa_{-n,n}(a).
}

This completely determines the operators $\fG_k(a)$ if $K=0$.
Our goal is to find a closed formula for all $g_\ep$,
which will imply a closed formula for the operators
$\fG_k(a)$ and $\fch_k(E)$ without the requirement that $K=0$.


\subsection{Lehn's boundary operator}
To control the intersection theory of Hilbert schemes,
in addition to the Fock space description of $\bH$,
we require Lehn's boundary operator,
which is closely related to the Laplace--Beltrami and Nabla operators
appearing in other contexts.
This operator admits an explicit description in terms of the
Heisenberg operators and the canonical class $K\in H$ of the surface.

Lehn's boundary operator \cite[Def.~3.8]{lehn_chern}
is defined by
\eq{\label{lehn d}
\ld:=\fch_1(\cO_X)=\fG_1(1_X)\in\iend_0(\bH),}
where $1=1_X\in H^0(X,\bQ)$.
We have $\ld(u)=\ch_1(\hb\cO n_X)\cdot u$ for $u\in \bH_n$.
For $f\in\End(\bH)$, define
\eq{f'=[\ld,f],\qquad
f^{(0)}=f,\qquad f^{(k+1)}=(f^{(k)})'.
}

\begin{theorem}[See {\lehn[Theorem 3.10]}]
\label{lehn1}
For $a\in H$ and $n\in\bZ$, we have
\eq{\label{lehn1:eq}
\fa'_n(a)=[\ld,\fa_n(a)]
=n\fL_n(a)-\frac{n(\n n-1)}2\fa_n(Ka),}
where $\fL_n(a)=\oh\sum_{i+j=n}\fa_{i,j}(a)$ \cref{oper L}.
\end{theorem}

\begin{remark}
Compared to \lehn,
we use operators $\fa_n=\fq^L_{-n}$ and $\fL_n=-\fL^{L}_{-n}$.
\end{remark}

\begin{corollary}
We have
$T=\fL_{-1}(1)=-[\ld,\fa_{-1}(1)]$.
\end{corollary}

\begin{corollary}[\cf {\lehn[Corollary 3.11]}]
\label{lehn:uniq}
Let $f\in\End(\bH)$ satisfy
\[
f\one=0,\qquad
[\ld,f]=[\fa_{-1}(a),f]=0\qquad\forall a\in H.\]
Then $f=0$.
\end{corollary}
\begin{proof}
Since $T=-[\ld,\fa_{-1}(1)]$, we have $[T,f]=0$.
Since $\fa_{-n-1}(a)=\frac1{n!}\ad_T^n\fa_{-1}(a)$,
we have $[\fa_{-n-1}(a),f]=0$ for all $n\ge0$ and $a\in H$.
Since $\bH$ is spanned by the elements of the form
$\fa_{-n_1-1}(a_1)\dots\fa_{-n_k-1}(a_k)\one$ for $n_i\ge0$ and $a_i\in H$, we obtain $f=0$.
\end{proof}


\begin{corollary}\label{ld expl:cor}
We have
\eq{\label{d-expl}
\ld=
-\frac16\sum_{\ov{\bn\in\Zs^3}{\n\bn=0}}\fa_\bn(1)
+\frac12\sum_{n>0}(n-1)\fa_{-n,n}(K).
}
\end{corollary}

\begin{proof}
Let $D$ be the right-hand side.
We have $\ld\one=0$ by degree reasons
and we have $D\one=0$, since $D$ is a combination
of normally ordered products.
Therefore to prove that $\ld=D$,
it suffices to show that
$[\ld-D,\fa_n(a)]=0$ for all $n<0$ and $a\in H$.
We have
\[
\frac16\sum_{\ov{\bn\in\Zs^3}{\n\bn=0}}[\fa_\bn(1),\fa_n(a)]
=[\W^3_0(1),\fa_n(a)]
=-n \W^2_n(a)
=-n\fL_n(a).
\]
For the quadratic term,
\[
\frac12\sum_{m>0}(m-1) \sbr{\fa_{-m,m}(K),\fa_n(a)}
=-\frac{n(\n n-1)}2\fa_n(Ka).
\]
Therefore
$
[D,\fa_n(a)]
=n\fL_n(a)-\frac{n(\n n-1)}2\fa_n(Ka)
=[\ld,\fa_n(a)]$.
\end{proof}

\begin{remark}
Let $\fa(z)=\sum_{n\ne0}\fa_nz^{-n-1}:H\to\cF(V)$.
Then \cref{d-expl} may be written as
\[
\ld=
[z^{-3}]\rbr{
-\frac1{6}\nop{\fa(z)^3}(1)
+\frac12\rbig{\dd_z\fa(z)_-\,\fa(z)_+}(K)
}.
\]
\end{remark}

The tautological operator $\fG_k(a)$ for $k\ge0$
satisfies
\eq{\label{d-G-comm}
\fG_k(a)\one=0,\qquad [\ld,\fG_k(a)]=0.}
In view of \cref{lehn:uniq},
it is uniquely determined by the commutators
$[\fG_k(a),\fa_{-1}(b)]$.
The following result is proved in \lehn[Theorem 4.2]
for Chern character operators
and in \cite[Lemma 5.8]{li_vertex} for tautological operators.

\begin{theorem}\label{lehn2}
For $a,b\in H$ and $k\ge0$, we have
\eq{\label{lehn2:eq}
[\fG_k(a),\fa_{-1}(b)]=\frac1{k!}\pr\fa k_{-1}(ab).
}
\end{theorem}

In the next section, we will interpret
\cref{d-G-comm,lehn2:eq} as a system of differential equations.
By solving this system, we will determine the operators $\fG_k(a)$.

\section{Differential equations for tautological operators}
\label{diff:sec}
Our goal is to obtain explicit formulas for the universal coefficients of tautological operators.
The key observation is that their commutation relations can be expressed in terms of differential operators on a polynomial algebra. This leads to a system of differential equations whose solution determines the desired coefficients.
As before, we consider the Frobenius algebra $H=H^*(X,\bQ)$
and the Fock space $V=V(H)$ isomorphic to $\bH=\bop_{n\ge0}H^*(\hb Xn,\bQ)$.

\subsection{Differential boundary operator}
Recall that for $a\in H$, we define \cref{N-A}
\[\N_a:\cA_H\to\End(\bH),\qquad
bx_\la=x_\la\ts b\mto \fa_\la(ba).\]
We use Lehn's commutator formula \cref{lehn1:eq} for $[\ld,\fa_n(a)]$
to describe the commutator $[\ld,\N_a(f)]$ in terms of differential operators on $\cA_H$.
Define the differential operators
\begin{equation}\label{eq:operators}
\begin{gathered}
\cM=\frac12\sum_{i,j\in\Z}(i+j)x_ix_j\dd_{i+j},\qquad
\cW=\frac12\sum_{i\in\Z}i(\n i-1)x_i\dd_i,\\
\cC=\oh\sum_{ij>0}\n ij\, x_{i+j}\dd_i\dd_j.
\end{gathered}
\end{equation}
We call them \emph{merge, weight and contraction operators}, respectively.

\begin{theorem}\label{diff bound}
For $a\in H$ and $f\in\cA_H$, we have
\eq{\label{D}
[\ld,\N_a(f)]=\N_a(\cD f),\qquad
\cD=\cM-K\cW+e\cC.
}
\end{theorem}

\begin{proof}
It suffices to consider $f=x_\la$, so that $\N_a(f)=\fa_\la(a)$.
Write $\la=\mset{n_1,\ldots,n_r}$, where $n_1\le\cdots\le n_r$.
By \cref{lehn1}, $[\ld,\fa_\la(a)]$ is obtained
by replacing each factor $\fa_{n_s}(u)$ by
\[
\frac{n_s}{2}\sum_{k+l=n_s}\nop{\fa_k\fa_l}(u)
-\frac{n_s(|n_s|-1)}2\fa_{n_s}(Ku).
\]
Note that if $\de(a)=u\ts v$, then
$\de(Ka)=Ku\ts v=u\ts Kv$ by \cref{frob:eq}.

The resulting products need not be normally ordered.
We restore normal order using the Heisenberg relations.
The terms requiring no additional interchanges give
\[
\N_a\bigl((\cM-K\cW)x_\la\bigr).
\]

Additional contractions arise when a newly introduced mode
must be commuted past an original mode with the opposite index.
Since the original modes are normally ordered, this can happen
only when a new negative mode crosses an earlier positive mode,
or a new positive mode crosses a later negative mode.
In either case, the two original modes have the same sign.

More precisely, for $s<t$ with $n_s,n_t>0$, replacing $n_t$
by $(-n_s,n_s+n_t)$ produces a contraction with coefficient
$n_sn_t$. Similarly, for $n_s,n_t<0$, replacing $n_s$ by
$(n_s+n_t,-n_t)$ produces a contraction with coefficient
$-n_sn_t$.
By \cref{frob contr}, contracting
the corresponding two factors of $\de_{r+1}(a)$ gives
$\de_{r-1}(ea)$. Thus the total contraction contribution is
\[
\sum_{\substack{s<t\\n_sn_t>0}}
\sgn(n_s)n_sn_t\,
\fa_{\la\ms\set{n_s,n_t}\sqcup\set{n_s+n_t}}(ea)
=\N_a(e\cC x_\la).
\]
Combining the contributions proves the assertion.
\end{proof}

Recall the translation operator \cref{trans}
\[\cT=-\sum_{n\in\bZ}nx_{n-1}\dd_n.
\]

\begin{lemma}
We have
\[[\cM,\dd_1]=\cT,\qquad
[\cW,\dd_1]=[\cC,\dd_1]=0.
\]
Therefore $[\cD,\dd_1]=\cT$.
\end{lemma}
\begin{proof}
We have
$[\dd_1,\cM]=\oh\sum_{i,j}(i+j)\dd_1(x_ix_j)\dd_{i+j}
=\sum_i(i+1)x_{i}\dd_{i+1}=-\cT$.
Similarly,
$[\dd_1,\cW]=\oh\sum_i i(\n i-1)\dd_1(x_i)\dd_i=0$
and $[\dd_1,\cC]=\oh\sum_{ij>0}\n i j\dd_1(x_{i+j})\dd_i\dd_j=0$.
\end{proof}

\begin{lemma}
Let $f\in\bop_{r>0}\cA_{r,0}$ satisfy $\cD f=\dd_1f=0$.
Then $f=0$.
\end{lemma}
\begin{proof}
We have $[\cD,\dd_1]=[\cM,\dd_1]=\cT$ and $[\cT,\dd_n]=(n+1)\dd_{n+1}$ for $n\ge1$.
Therefore $\dd_n f=0$ for $n\ge1$,
hence $f$ depends only on $x_{-1},x_{-2},\dots$.
Since $f$ has charge zero, we obtain $f=0$.
\end{proof}

\begin{remark}
The above result is an analogue of \cref{lehn:uniq}.
\end{remark}

\subsection{A system of differential equations}
For $f\ts \ep\in\cA_{r,n}\ts H$,
define its total degree $r+\n \ep/2$ and charge~$n$.
Let $\cE=\set{1,e,K,K^2}$ and let
$\Ga(t)
=\sum_{k\ge0}\Ga_kt^k$,
where
\eq{
\Ga_k=\sum_{\ep\in\cE}\ep g_{\ep,k}\in\cA_H
}
has total degree $k+2$ and charge $0$, and satisfies
(\cf \cref{LQW:general})
\eq{
\fG_k(a)=(-1)^k\N_a(\Ga_k).
}
Then
\[
\N_a(\cD\Ga_k)=(-1)^k[\ld,\fG_k(a)]=0.
\]
Applying \cref{N-comm,lehn2}, we obtain
\[
\N_{ab}(\dd_1\Ga(t))=
\sum_{k\ge0}[\fG_k(a),\fa_{-1}(b)](-t)^k
=e^{-t\ad_\ld}(\fa_{-1}(ab))
=\N_{ab}(e^{-t\cD}x_{-1}).
\]
Moreover
$\N_a(\Ga(0))=\fG_0(a)=\fL_0(a)=\N_a\rbr{\sum_{n>0}x_{-n}x_n}$.

Motivated by these identities, we seek a series $\Ga(t)$ satisfying the following system of differential equations
\eq{\label{system1}
\cD\Ga=0,\qquad
\dd_1\Ga=e^{-t\cD}x_{-1}, \qquad
\Ga(0)=\sum_{n>0}x_{-n}x_n.
}
It is equivalent to the system
\eq{\label{system2}
\cD\Ga=0,\qquad (\cT+\dd_1\dd_t)\Ga=0, \qquad
\Ga(0)=\sum_{n>0}x_{-n}x_n.
}
Indeed, if \cref{system1} holds,
then $[\cD,\dd_1]=\cT$ implies
$\cT\Ga=\cD\dd_1\Ga=-\dd_t\dd_1\Ga$.
If \cref{system2} holds,
then
$\dd_t(\dd_1\Ga)=-\cT\Ga=-[\cD,\dd_1]\Ga=-\cD(\dd_1\Ga)$.
Since $\Ga(0)=\sum_{n>0}x_{-n}x_n$, we have $\dd_1\Ga(0)=x_{-1}$.
This implies $\dd_1\Ga=e^{-t\cD}x_{-1}$.

\subsection{Solution strategy}
We use the following strategy to solve~\cref{system2}.
Let
\eq{
\Ga(t)=\sum_{\ep\in\cE}\ep \hg_\ep(t)
=\hg_1(t)+e \hg_e(t)+K \hg_K(t)+K^2 \hg_{K^2}(t).}
Since $K^3=Ke=e^2=0$,
to solve the equation $\cD\Ga=0$, where $\cD=\cM-K\cW+e\cC$,
it suffices to solve the following four equations:
\eq{
\cM \hg_1=0,\qquad
\cM \hg_e=-\cC \hg_1,\qquad
\cM \hg_K=\cW \hg_1,\qquad
\cM \hg_{K^2}=\cW \hg_K.
}
Additionally, we require $(\cT+\dd_1\dd_t)\hg_\ep=0$ for all $\ep\in\cE$ and $\Ga(0)=\sum_{n>0}x_{-n}x_n$.

We first construct $\hg_1$, then $\hg_e$ and $\hg_K$, and finally $\hg_{K^2}$.
At each step, we seek a particular solution of an equation
involving $\cM$, determined up to an element of $\Ker\cM$.
The translation equation and the initial condition are then used to fix this ambiguity.

Given a Laurent series $f(z)$, define
\eq{\ct f=[z^0]f,\qquad
\de=z\dd_z,\qquad
\te=\n\de:z^n\mto \n nz^n.
}

\begin{lemma}
We have
\begin{gather}
\ct{f\cdot \de g}=-\ct{\de f\cdot g},\qquad
\ct{f\cdot \te g}=\ct{\te f\cdot g},\label{de-prop}
\\
\ct{\te f\cdot\de f}=0,\qquad \ct{z\te f}=\ct{zf}.
\end{gather}
\end{lemma}
\begin{proof}
Let $f=\sum_i f_iz^i$ and $g=\sum_i g_iz^i$.
Then
\[
\ct{f\te g}
=\sum_{i+j=0}\n j f_ig_j
=\sum_{i+j=0}\n i f_ig_j=\ct{g\te f}.
\]
We also have $\ct{f\de g+g\de f}=\ct{\de(fg)}=0$.
Then
\[\ct{\te f\cdot\de f}
=\ct{f\cdot \te\de f}=-\ct{\te\de f\cdot f},
\]
hence $\ct{\te f\cdot\de f}=0$.
Finally, $\ct {z\te f}=\ct{\te(z)f}=\ct{zf}$.
\end{proof}

\begin{lemma}
We have
\eq{
\ct{z\rbig{(\te f)\de g+(\de f)\te g}}=0,\qquad
\ct{z(\te f)\de f}=0.
}
\end{lemma}
\begin{proof}
Let $f=\sum_if_iz^i$ and $g=\sum_i g_iz^i$.
Then
\[\ct{z\rbig{(\te f)\de g+(\de f)\te g}}
=\sum_{i+j=-1}(\n ij+i\n j)f_ig_j=0,
\]
since $i+j=-1$ implies $ij\le0$, hence $\n ij+i\n j=0$.
Taking $f=g$, we obtain $\ct{z(\te f)\de f}=0$.
\end{proof}

\begin{lemma}
We have
\[\ct{(\te f)^2\de f}
=-\frac13\ct{(\de f)^3}.\]
\end{lemma}
\begin{proof}
Let $f=f_++f_-$, where $f_+$ has only nonnegative powers of $z$ and $f_-$ has only negative powers of $z$.
Using $a=\de f_+$ and $b=\de f_-$, we have $\te f=a-b$.
Therefore
\[
(\te f)^2\de f+\frac13(\de f)^3
=\frac13(4a^2+4b^2-4ab)(a+b)=\frac43(a^3+b^3).
\]
But $\ct{a^3}=\ct{b^3}=0$.
\end{proof}

\subsection{First identities}
We begin by establishing several differential identities for exponential generating functions.
These identities will provide the building blocks for the solution of our system of differential equations.
Let $X=\sum_{n\ne0}x_nz^n$ and let
\begin{gather}
V=\ct{e^{tX}}-1,\qquad
P=t\ct{e^{tX}\de^2X},\qquad Y=\ct{ze^{tX}},\\
R_e=\ct{e^{tX}\de^3X},\qquad
R_K=\ct{e^{tX}\de\te X}.
\end{gather}
We have
\begin{gather}
\cM X=\oh\sum_{i,j}(i+j)x_ix_j z^{i+j}
=\oh\de(X^2)=X\cdot \de X,\\
\cW X=\frac12\sum_{n} n(\n n-1)x_nz^n=\oh\de(\te-1)X.
\label{WX}
\end{gather}

\begin{lemma}
We have
\begin{gather}
\cM V=0,\qquad \cM P=-2R_e,\\
\cC V=\frac{t^2}{12}R_e,\qquad \cW V=\frac t2R_K,\\
(\cT+\dd_1\dd_t)V=2Y,\qquad
(\cT+\dd_1\dd_t)P=4Y.
\end{gather}
\end{lemma}

\begin{proof}
\begin{clm-l}
\item
Since $\cM X=X\de X$, we have
\[\cM V=\ct{\cM(e^{tX})}
=t\ct{Xe^{tX}\de X}
=\ct{\de\rbr{Xe^{tX}-e^{tX}/t}}=0.
\]

\item
For $P=t\ct{e^{tX}\de^2X}$, we have
\begin{multline*}
\cM P
=t\ct{\cM(e^{tX}\de^2X)}
=t\ct{e^{tX}\left(tX(\de X)(\de^2X)
+\de^2(X\de X)\right)}\\
=t\ct{\de(Xe^{tX}\de^2X)
+2e^{tX}(\de X)(\de^2X)}\\
=2t\ct{e^{tX}(\de X)(\de^2X)}
=2\ct{\de(e^{tX}\de^2 X)-e^{tX}\de^3X}
=-2R_e.
\end{multline*}

\item
We have $\cC=\frac12\sum_{ij>0}\n i jx_{i+j}\dd_i\dd_j$.
Since $\dd_i\dd_j e^{tX}=t^2z^{i+j}e^{tX}$, we obtain
\[\cC e^{tX}
=\frac{t^2}2e^{tX}\sum_{ij>0}\n i jx_{i+j}z^{i+j}.
\]
We have $\sum_{{i+j=n},{ij>0}}\n i j=\frac16(n^3-n)$,
hence
\[\cC e^{tX}
=\frac{t^2}{12}e^{tX}(\de^3-\de)X.
\]
We have $\ct{e^{tX}\de X}=\frac1t\ct{\de e^{tX}}=0$.
Therefore
\[\cC V=\ct{\cC e^{tX}}=\frac{t^2}{12}\ct{e^{tX}\de^3X}=\frac{t^2}{12}R_e.\]

\item
We have $\cW=\frac12\sum_{i\in\Z}i(\n i-1)x_i\dd_i$.
Therefore $\cW X=\oh(\de\te-\de)X$ and
\[\cW V=\ct{\cW e^{tX}}
=\frac t2 \ct{e^{tX}(\de\te-\de)X}
=\frac t2\ct{e^{tX}\de\te X}=\frac t2 R_K.
\]

\item
We have $\cT=-\sum_{n\in\bZ}nx_{n-1}\dd_n$.
Therefore $\cT X=-z(1+\de)X$ and
\[\cT V=-t\ct{ze^{tX}(X+\de X)}\]
Since $\dd_1X=z$, we obtain
$\dd_1\dd_t V=\ct{\dd_1 (Xe^{tX})}=\ct{z(1+tX)e^{tX}}$.
Therefore
\[(\cT+\dd_1\dd_t)V
=\ct{ze^{tX}(1-t\de X)}
=\ct{2ze^{tX}-\de(ze^{tX})}
=2\ct{ze^{tX}}=2Y.
\]

\item
We have $\cT X=-z(1+\de)X$ and $\de (zf)=z(1+\de)f$, hence
\[\cT(\de^2X)
=\de^2(\cT X)=-z(1+\de)^3X.
\]
For $P=t\ct{e^{tX}\de^2X}$, we obtain
\[
\begin{aligned}
\cT P
&=-t\ct{ze^{tX}
\left(t(X+\de X)\de^2X+(1+\de)^3X\right)},\\
\dd_1\dd_tP
&=\ct{ze^{tX}
\left(1+tX+2t\de^2X+t^2X\de^2X\right)}.
\end{aligned}
\]
Therefore
\begin{multline*}
(\cT+\dd_1\dd_t)P
=\ct{ze^{tX}(1-3t\de X)
-t\de(ze^{tX}\de^2X)}\\
=\ct{ze^{tX}(1-3t\de X)}
=\ct{4ze^{tX}-3\de(ze^{tX})}
=4Y.
\tag*{\qedhere}
\end{multline*}
\end{clm-l}
\end{proof}

\subsection{Logarithmic identities}
The contributions involving the canonical class require additional generating functions, naturally expressed through logarithmic series.
We derive the identities needed to construct these contributions and complete the solution of the differential equations.
Let $X=\sum_{n\ne0}x_nz^n$ and let $\Om=-\log(1-uX)$.
Define
\begin{gather}
\hB(u)=\ct{\Om\te\Om}=\sum_{r\ge2}B_ru^r,\qquad
B(t)=\dd_t\cB(\hB)
=\cB(u\i\hB)
=\sum_{r\ge2}\frac{B_r}{(r-1)!}t^{r-1},\\
\hD(u)=\ct{e^\Om(\te\Om)^2}=\sum_{r\ge2}D_ru^r,\qquad
D(t)
=\cB(\hD)
=\sum_{r\ge2}\frac{D_r}{r!}t^{r},
\end{gather}
where $\cB(u^r)=t^r/r!$ is the Borel transform.

\begin{lemma}
We have
\begin{gather}
\cM B=-\frac2t R_K,\qquad \cM D=-2\cW B+\frac13R_e,\\
(\cT+\dd_1\dd_t)B=\frac2tY,
\qquad
(\cT+\dd_1\dd_t)D=\dd_1B.
\end{gather}
\end{lemma}

\begin{proof}
We have $e^{-\Om}=1-uX$ and $\cM X=X\de X$, hence
\[
\cM\Om=e^\Om u \cM X
=e^\Om (1-e^{-\Om})\de X
=u\i\de\Om-\de X.
\]
Let $\cS=\cT+u\i\dd_1$.
We have $\cT X=-z(1+\de)X$ and $\dd_1X=z$, hence
\[\cS\Om
=e^\Om(u\cT X+\dd_1X)=ze^\Om\rbr{1-uX-u\de X}
=z(1-\de\Om).
\]
We have $\hB=\ct{\Phi(\Om,\te\Om)}$ for $\Phi(x,y)=xy$
and $\hD=\ct{\Phi(\Om,\te\Om)}$ for $\Phi(x,y)=e^xy^2$.
Let $\Phi_x=\dd_x\Phi$ and $\Phi_y=\dd_y\Phi$.
For a derivation $\bde$ of $\cA$, we have
$\bde\ct{\Phi(\Om,\te\Om)}
=\ct{\rbig{\Phi_x+\te\Phi_y}\bde\Om}$.
Hence
\[\cM\ct{\Phi}
=\ct{u\i\de\Phi-(\Phi_x+\te\Phi_y)\de X}
=-\ct{(\Phi_x+\te\Phi_y)\de X},\]
\begin{multline*}
\cS\ct{\Phi}
=\ct{z(1-\de\Om)(\Phi_x+\te\Phi_y)}
=\ct{z\rbig{\Phi_x+\te\Phi_y-\Phi_x\de\Om+(\de\Phi_y)\te\Om}}
\\
=\ct{z\rbig{\Phi_x+\te\Phi_y+\de(\Phi_y\te\Om)-\de\Phi}}
=\ct{z\rbig{\Phi_x+\Phi_y-\Phi_y\te\Om+\Phi}}.
\end{multline*}
This implies
\begin{gather*}
\cM\hB
=-2\ct{(\te\Om)(\de X)}.
\\
\cS\hB
=\ct{z\rbr{\te\Om+\Om-\Om\te\Om+\Om\te\Om}}=2\ct{z\Om}.
\\
\cM\hD
=-\ct{\rbr{e^\Om(\te\Om)^2+2\te(e^\Om\te\Om)}\de X}
=-\ct{e^\Om\rbr{(\te\Om)^2 \de X+2(\te\Om)\de\te X}}.
\\
\cS\hD
=\ct{ze^\Om\rbr{(\te\Om)^2+2\te\Om-2(\te\Om)^2+(\te\Om)^2}}
=2\ct{ze^\Om\te\Om}.
\end{gather*}

By \cref{WX}, we have $\cW\Om=e^\Om u\cW X=\frac 12e^\Om\de(\te-1)uX$,
hence
\[\cW \hB
=2\ct{(\te\Om)\cW\Om}
=u\ct{e^\Om(\te\Om)\de(\te-1)X}.
\]
Therefore
\begin{multline*}
\cM\hD+\frac2u\cW\hB
=-\ct{e^\Om\rbr{(\te\Om)^2+2\te\Om}\de X}
=-\frac1u\ct{\rbr{(\te\Om)^2+2\te\Om}\de\Om}\\
=-\frac1u\ct{(\te\Om)^2\de\Om}
=\frac1{3u}\ct{(\de\Om)^3}.
\end{multline*}
We have $u\de^3 X=-\de^3(e^{-\Om})=e^{-\Om}
\rbr{(\de\Om)^3-3(\de\Om)(\de^2\Om)+\de^3\Om}$.
Since $\ct{(\de\Om)\de^2\Om}=0$ and $\ct{\de^3\Om}=0$,
we obtain
\[\cM\hD+\frac2u\cW\hB
=\frac1{3u}\ct{(\de\Om)^3}
=\frac1{3}\ct{e^\Om\de^3X}
=\frac1{3}\ct{\frac{\de^3X}{1-uX}}.
\]
We also have $\dd_1\hB=2\ct{(\te\Om)\dd_1\Om}=2u\ct{ze^\Om\te\Om}$, hence $\cS\hD=u\i\dd_1\hB$.

Since the Borel transform satisfies
$\cB(u\i\Om)=\frac{e^{tX}-1}{t}$ and
$\cB(e^\Om)=e^{tX}$,
we obtain
\begin{gather*}
\cM B
=\cB(u\i\cM\hB)
=-2\ct{\cB(u\i\Om)\te\de X}
=-\frac2tR_K,\\
(\cT+\dd_1\dd_t)B
=\cB(u\i\cS\hB)
=2\ct{z\cB(u\i\Om)}
=\frac2tY,\\
\cM D+2\cW B
=\cB\rbr{\cM\hD+\frac2u\cW\hB}
=\frac13\ct{\cB(e^\Om)\de^3X}
=\frac13R_e,\\
(\cT+\dd_1\dd_t)D
=\cB(\cS\hD)
=\cB(u\i\dd_1\hB)
=\dd_1B.
\tag*{\qedhere}
\end{gather*}
\end{proof}

\subsection{Closed formulas for tautological operators}
We now combine the identities established above to construct an explicit solution of the differential equations \cref{system2}.
This yields a closed formula for the tautological operators.

\begin{theorem}\label{sol:th1}
Consider the series $\Ga(t)=\sum_{k\ge0}\Ga_kt^k\in\cA_H\pser t$ given by
\begin{equation}
\Ga=\frac{V}{t^2}
+\frac e{24}\rbr{P-2V}
+K\rbr{\frac{V}{2t}-\frac B4}
+K^2\rbr{\frac{V}{12}+\frac{P}{48}-\frac{tB}8+\frac D8}.
\end{equation}
Then
\[\cD\Ga=0,\qquad (\cT+\dd_1\dd_t)\Ga=0,\qquad
\Ga(0)=\sum_{n>0}x_{-n}x_n.
\]
Moreover, $\fG_k(a)=(-1)^k\N_a(\Ga_k)$.
\end{theorem}
\begin{proof}
Let us write the summands of $\Ga$ as
$\Ga=\hg_1+e\hg_e+K\hg_K+K^2\hg_{K^2}$.
Then
\[\cD\Ga=\cM\hg_1
+e(\cM\hg_e+\cC\hg_1)
+K(\cM\hg_K-\cW\hg_1)
+K^2(\cM\hg_{K^2}-\cW\hg_K).
\]
We have $\cM\hg_1=\frac1{t^2}\cM V=0$,
\[
\cM\hg_e=\frac1{24}\cM P=-\frac1{12}R_e=-\cC\hg_1,
\qquad
\cM\hg_K=-\frac14\cM B=\frac1{2t}R_K=\cW\hg_1,\]
\begin{multline*}
\cM\hg_{K^2}
=\frac1{48}\cM P-\frac t8\cM B+\frac18\cM D\\
=-\frac1{24}R_e+\frac14R_K-\frac14\cW B+\frac1{24}R_e
=\frac1{2t}\cW V-\frac14\cW B=\cW\hg_K.
\end{multline*}
This implies that $\cD\Ga=0$.
Let $\cS=\cT+\dd_1\dd_t$ and $Y=\ct{ze^{tX}}$.
We have
\[\cS V=2Y,\qquad \cS P=4Y,\qquad \cS B=\frac{2Y}t,\qquad
\cS D=\dd_1B.
\]
For a scalar series $h(t)$, we have
$\cS(hV)
=h\cS V+h'\dd_1V$ and $\dd_1V=tY$.
Therefore
\[\cS\hg_1=t^{-2}\cS V-2t^{-3}\dd_1V
=2t^{-2}Y-2t^{-2}Y=0.
\]
\[\cS\hg_e=\frac1{24}(4Y-4Y)=0,\qquad
\cS\hg_K=\frac1{2t}\cS V-\frac1{2t^2}\dd_1V-\frac14\cS B
=\frac Yt\rbr{1-\frac12-\frac12}=0.
\]
\[\cS\hg_{K^2}
=\frac Y6+\frac Y{12}-\rbr{\frac Y4+\frac18\dd_1B}+\frac18\dd_1B=0.
\]
Therefore $\cS\Ga=0$.

We have $B=O(t)$, $D=O(t^2)$, $V=O(t^2)$, $P=O(t^2)$, hence
\[\Ga(0)=[t^2]V=\frac12\ct{X^2}=\sum_{n>0}x_{-n}x_n.\]
We have seen that the above equations imply $\dd_1\Ga=e^{-t\cD}x_{-1}$
(see \cref{system2}).

Let $\tG_k(a)=(-1)^k\N_a(\Ga_k)$.
It follows from \cref{diff bound} that \[[\ld,\tG_k(a)]=(-1)^k\N_a(\cD\Ga_k)=0.\]
It follows from \cref{N-comm} and \cref{lehn2} that
\[[\tG_k(a),\fa_{-1}(b)]
=(-1)^k\N_{ab}(\dd_1\Ga_k)
=\frac1{k!}\N_{ab}\rbr{\cD^k x_{-1}}
=\frac1{k!}\ad_\ld^k(\fa_{-1}(ab))
=[\fG_k(a),\fa_{-1}(b)].
\]
Every $\tG_k(a)$ annihilates the vacuum, since
$\Ga_k$ is a linear combination of normally ordered
monomials of charge zero and positive length.
We conclude from \cref{lehn:uniq} that $\fG_k(a)-\tG_k(a)=0$.
\end{proof}

\begin{corollary}
\label{cor:g_ep}
Let $\la=\mset{n_1,\dots,n_r}\in\Pg_{r,0}$ with $r\ge2$.
For $I\sbs[r]$, let $n_I=\sum_{i\in I}n_i$.
Then the universal coefficients in \cref{LQW:general} are given by
\begin{gather}
g_1(\la)=1,\\
g_e(\la)=\frac{p_2(\la)-2}{24},\\
g_K(\la)=\frac12-\frac{B_{r}(\la)}{4(r-1)!},\\
g_{K^2}(\la)=\frac1{12}+\frac{p_2(\la)}{48}
-\frac{B_r(\la)}{8(r-1)!}+\frac{D_r(\la)}{8r!},
\end{gather}
where
\begin{gather}
B_r(\la)=\sum_{\ov{[r]=I\sqcup J}{I,J\ne\es}}
(\n {I}-1)!(\n J-1)!\n{n_I},\\
D_r(\la)=\sum_{\ov{[r]=I\sqcup J\sqcup J'}{I,J\ne\es}}
(\n {I}-1)!(\n J-1)!\n {J'}!
\n{n_I}\n{n_J}.
\end{gather}
Here the block $J'$ may be empty.
\end{corollary}

\begin{proof}
Let $x_\la=\prod_{i=1}^r x_{n_i}=\prod_{k\ne0}x_k^{m_k(\la)}$.
By the definition of the universal coefficients,
\[
\la![x_\la]\Ga(t)
=g_1(\la)t^{r-2}
+e\,g_e(\la)t^r
+K\,g_K(\la)t^{r-1}
+K^2g_{K^2}(\la)t^r.
\]
Thus it suffices to compute $\la![x_\la]$ on
$V,P,B,D$.
Consider the algebra morphism
\[\rho_\la:\bQ[x_n\col n\in\Zs]\to\bQ[y_1,\dots,y_r],\qquad
x_n\mto \sum_{\ov{i\in[r]}{n_i=n}}y_i.
\]
Writing $y_I=\prod_{i\in I}y_i$ for $I\sbs [r]$, we have
\[
\la![x_\la]f=[y_{[r]}]\rho_\la(f).
\]
Let $X=\sum_{n\ne0}x_nz^n$ and let
\[
X_\la=\rho_\la(X)=\sum_{i=1}^r y_iz^{n_i}.
\]
Since $n_{[r]}=0$, the definitions of $V=\ct{e^{tX}}-1$ and
$P=t\ct{e^{tX}\de^2X}$ give
\[
\la![x_\la]V
=[y_{[r]}]\ct{e^{tX_\la}}=t^r,\qquad
\la![x_\la]P
=[y_{[r]}]\ct{te^{tX_\la}\de^2X_\la}
=p_2(\la)t^r.
\]

Next, using
\[
\Om=-\log(1-uX)=\sum_{k\ge1}\frac{u^kX^k}{k},
\qquad
e^\Om=\sum_{k\ge0}u^kX^k,
\]
we obtain, for $I\ne\es$,
\[
\begin{aligned}
[y_I]\rho_\la(\Om)
&=(\n I-1)!u^{\n I}z^{n_I},\\
[y_I]\rho_\la(\te\Om)
&=(\n I-1)!\n{n_I}u^{\n I}z^{n_I},
\end{aligned}
\]
and, for arbitrary $I\sbs[r]$,
\[
[y_I]\rho_\la(e^\Om)
=\n I!u^{\n I}z^{n_I}.
\]
Consequently, the definitions of $\hB$ and $\hD$ yield
\[
\la![x_\la]\hB=u^rB_r(\la),\qquad
\la![x_\la]\hD=u^rD_r(\la).
\]

Applying the Borel transform, we obtain
\[
\la![x_\la]B
=\frac{B_r(\la)}{(r-1)!}t^{r-1},
\qquad
\la![x_\la]D
=\frac{D_r(\la)}{r!}t^r.
\]

Substituting these four identities into the formula
for $\Ga$ in \cref{sol:th1} and comparing with
the defining expansion of $\la![x_\la]\Ga(t)$ gives
the statement.
\end{proof}

\section{Quasi-modular forms}
\label{QM:sec}
We introduce the analytic tools needed to study the arithmetic properties of the generating series arising from vertex operators.
In particular, we develop the Fourier theory of quasi-elliptic functions and its relation to quasi-modular forms and graph sums.
The resulting Fourier coefficient formulas will play a central role in our proof of the main theorem.

\subsection{Quasi-modular forms}
Let $Q_0(q)=1$, $D_q=q\dd_q$ and
\eq{Q_k(q)=D_q^{k-1}\frac q{1-q}
=\sum_{n\ge1}n^{k-1}q^n,\qquad k\ge1.
}
Define the Eisenstein series (for even $k\ge1$)
\eq{
E_k(q)=1-\frac{2k}{B_k}\sum_{n,d\ge1}n^{k-1}q^{nd}
=1-\frac{2k}{B_k}F_k(q),
}
where $B_k\in\bQ$ are Bernoulli numbers
satisfying $\sum_{k\ge0}B_k\frac{t^k}{k!}=\frac{t}{e^t-1}$ and
\eq{F_k(q)
=\sum_{n\ge1}n^{k-1}\frac{q^n}{1-q^n}
=\sum_{n,d\ge1}n^{k-1}q^{nd}=\sum_{d\ge1}Q_k(q^d),
\qquad k\ge1.
}
The algebra of quasi-modular forms
\[\QM=\QM(1;\bQ)=\bQ[E_2,E_4,E_6]\]
contains all $E_{2k}$, hence all $F_{2k}$, for $k\ge1$.

\begin{remark}
Series of the form $\sum_{n\ge1}a_n\frac{q^n}{1-q^n}$ are called Lambert series.
We will encounter many such series below.
\end{remark}

\subsection{Quasi-elliptic algebras}
We refer to \cite{hongminh_quasia} for more details on the material below.
Define the $q$-Pochhammer symbol
\[
\ph z=\prod_{i\ge0}(1-q^iz)=\exp\rbr{-\sum_{n>0}\frac1n\frac{z^n}{1-q^n}}
\]
and define the theta function
\[
\Te(z,q)=\frac{\ph z\ph{z\i q}}{\ph q^2}
=\exp\rbr{\sum_{n>0}\frac1n\frac{2q^n-z^n-z^{-n}q^n}{1-q^n}}
\]
which has zeros only at $z=q^n$ for $n\in\bZ$.
Let $D_z=z\dd_z$ and let (\cf \cite[(3.9)]{zhu_modular})
\begin{equation}
\Phi_k(z,q)
=-D_z^k\log\Te(z,q)
=\sum_{n>0}\frac{n^{k-1}z^n-(-n)^{k-1}z^{-n}q^n}{1-q^n},\qquad k\ge1.
\end{equation}
%
%
The above expansion is valid on the annulus $\n q<\n z<1$.
We have
\begin{gather}
\Phi_1(qz,q)=\Phi_1(z,q)+1=-\Phi_1(z\i,q),\\
\Phi_k(qz,q)=\Phi_k(z,q)=(-1)^k\Phi_k(z\i,q),\qquad k\ge2,\\
\ra_k(n)=\what\Phi_k(n):=[z^n]\Phi_k(z,q)
=\frac{n^{k-1}}{1-q^n},\qquad k\ge1,\, n\ne0.\label{hPhi}
\end{gather}
%
We also use $\rho(n)=\rho_1(n)=\frac1{1-q^n}$ satisfying $\rho(n)+\rho(-n)=1$.
We have \cite{hongminh_quasia}
\eq{\label{Phi Laurent}
\Phi_k(e^s,q)\in\QM\lser s,\qquad k\ge1.}
Let $\cC_0=\QM$.
For $r\ge1$, put $z_0=1$ and define the quasi-elliptic algebra
\eq{
\cC_r
=\cC_0\sbr{\Phi_k(z_j/z_i)\col 0\le i\ne j\le r,\, k\ge1}\\
=\cC_0\sbr{\Phi_k(z_i),\Phi_k(z_j/z_i)\col 1\le i\ne j\le r,\, k\ge1}.
}
The elements of $\cC_r$ are meromorphic on
\[X_r=(\bC^\xx)^r\xx\sets{q\in\cC}{0<\n q<1}\]
and holomorphic on
\[
U_r=\sets{(z_1,\dots,z_r,q)\in\bC^{r+1}}
{1>\n{z_1}>\dots>\n{z_r}>\n q>0}\sbs X_r.
\]

We have $\Phi_k(z,q)\in\bQ(z)\pser q$ and
\[\cC_r\sbs\bQ(z_1,\dots,z_r)\pser q.\]
Consider the expansion map
\[
\io_r:\cC_r\emb\bQ\lser {z_1}\dots\lser {z_r}\pser q
\]
corresponding to $1>\n {z_1}>\dots>\n {z_r}>\n q>0$.

\begin{theorem}[Constant term theorem,
see {\cite{hongminh_quasia,goujard_counting}}]
\label{CTT}
For every $f\in\cC_r$, we have
\eq{
\CT_r(f):=[z_1^0\dots z_r^0]\io_r(f)\in\cC_0.
}
Moreover, $[z_r^0]\io_r(f)\in\io_{r-1}\cC_{r-1}$.
\end{theorem}

\begin{lemma}\label{Phi:local}
Let $f\in\cC_r$ and let $0\le i<r$.
Then
\[
f(z_1,\dots,z_{r-1},z_ie^s)\in\cC_{r-1}\lser s.
\]
\end{lemma}

\begin{proof}
It suffices to consider the generators $\Phi_k(z_b/z_a)$
of $\cC_r$, where $0\le a<b\le r$.
Generators not involving $z_r$ belong to $\cC_{r-1}$.
If a generator involves both $z_r$ and $z_i$,
then after the substitution $z_r=z_ie^s$,
we obtain $\Phi_k(e^s)\in\cC_0\lser s$.
Every other generator involving $z_r$
produces $\Phi_k(w e^{s})$,
where $w=z_i/z_j$ for some $0\le j<r$ with $j\ne i$.
Its Taylor coefficients belong to $\cC_{r-1}$, since
\[
\Phi_k(we^{s})
=\sum_{\ell\ge0}\frac{s^\ell}{\ell!}\Phi_{k+\ell}(w)
\in\cC_{r-1}\pser s.
\qedhere
\]
\end{proof}

\subsection{Graph sums}
We now relate the theory of quasi-elliptic functions to graph sums. The main observation is that conservation conditions on integer flows can be expressed through constant terms of products of quasi-elliptic functions, providing a natural explanation for the quasi-modularity of certain graph sums.
As a byproduct, we generalize some of the results in~\GM.

A quiver (or digraph)
is a tuple $Q=(Q_0,Q_1,s,t)$,
where $Q_0,Q_1$ are sets,
called sets of vertices and arrows,
and $s,t:Q_1\to Q_0$ are maps, called source and target maps.
We write $a:i\to j$ for $a\in Q_1$ with $s(a)=i$ and $t(a)=j$.
A path $p$ in $Q$ is a sequence of vertices and arrows
$i_0\xto{a_1}i_1\to\dots\to i_{n-1}\xto{a_n}i_n$ in $Q$.
We
define $s(p)=i_0$, $t(p)=i_n$,
and the length $\ell(p)=n$.
We say that $p$ is a cycle if $\ell(p)>0$ and $i_0=i_n$.
A quiver is called acyclic if it has no cycles.
A quiver may be interpreted as a $1$-dimensional semi-simplicial
set with the face maps $d_0=t$ and $d_1=s$.

Assume that $Q$ is finite, meaning that $Q_0,Q_1$ are finite.
Define the boundary operator (or divergence)
\eq{\label{boundary}
\dd:\bZ^{Q_1}\to\bZ^{Q_0},\qquad
\dd_i(\bn)=\sum_{t(a)=i}n_a-\sum_{s(a)=i}n_a.
}
We may interpret $\bn\in\bZ^{Q_1}$ as an integer flow on $Q$,
with $\dd_i(\bn)$ measuring the net inflow at vertex $i\in Q_0$.
Thus $\dd(\bn)=0$ expresses the conservation law at every vertex.

\begin{theorem}
\label{QM acyclic}
Let $Q$ be an acyclic quiver
and let $P\in\bQ[n_a\col a\in Q_1]$ and $f_a\in\cC_1$ for $a\in Q_1$.
Then
\[\sum_{\ov{\bn\in\Zs^{Q_1}}{\dd(\bn)=0}}
P(\bn)
\prod_{a\in Q_1} \what f_a(n_a)\in\QM,\]
where $\what f(n)=[z^n]f(z)$ for $f\in\cC_1$.
\end{theorem}
\begin{proof}
We can assume that $P$ is a monomial $\prod_{a\in Q_1} n_a^{k_a}$.
Then $n_a^{k_a}\what f_a(n_a)
=[z^{n_a}](D_z^{k_a}f_a)$,
hence we can substitute $f_a$ by $D_z^{k_a}f_a\in\cC_1$ for all $a\in Q_1$ and assume that $P=1$.
We have $c_a=\what f_a(0)\in\cC_0=\QM$,
hence we can substitute $f_a$ by $f_a-c_a$ and assume that $\what f_a(0)=0$.

Consider the variables $z_i$ for $i\in Q_0$ and let $z_a=z_j/z_i$ for $a:i\to j$.
We order $Q_0$ so that $i<j$ for every arrow $a:i\to j$.
Then $\prod_{a\in Q_1} f_a(z_a)
=\sum_{\bn\in\bZ^{Q_1}}\prod_{a\in Q_1}\what f_a(n_a)z_a^{n_a}$, hence
\[
\sum_{\ov{\bn\in\bZ^{Q_1}}{\dd(\bn)=0}}
\prod_{a\in Q_1} \what f_a(n_a)
=\CT\rbr{\prods_{a\in Q_1} f_a(z_a)}\in\cC_0
\]
by the constant term theorem.
\end{proof}

\begin{corollary}
\label{QM for rho}
Let $Q$ be an acyclic quiver
and let $P\in\bQ[n_a\col a\in Q_1]$.
Then
\[\sum_{\ov{\bn\in\Zs^{Q_1}}{\dd(\bn)=0}}
P(\bn)
\prod_{a\in Q_1} \frac1{1-q^{n_a}}\in\QM.\]
\end{corollary}
\begin{proof}
We apply \cref{QM acyclic} to $f_a=\Phi_1$ for $a\in Q_1$.
They satisfy $\what f_a(n)=\frac1{1-q^n}$ for $n\ne0$.
\end{proof}

\begin{corollary}
\label{QM for E_-}
Let $Q$ be a quiver with a decomposition
$Q_1=E_+\sqcup E_-$ such that reversing the arrows from $E_-$ gives an acyclic quiver.
Let $P\in\bQ[n_a\col a\in Q_1]$.
Then
\[\sum_{\ov{\bn\in\Zs^{Q_1}}{\dd(\bn)=0}}
P(\bn)\prod_{a\in E_+}\rho(n_a)\prod_{a\in E_-}\rho(-n_a)\in\QM,
\]
where $\rho(n)=\frac1{1-q^n}$ for $n\ne0$.
\end{corollary}
\begin{proof}
Let $Q'$ be the acyclic quiver obtained from $Q$ by reversing the arrows from $E_-$.
For $\bn\in\Zs^{Q_1}$, let $\bm\in\Zs^{Q'_1}$ be defined by $m_a=n_a$ for $a\in E_+$ and $m_a=-n_a$ for $a\in E_-$.
Then $\dd_Q(\bn)=\dd_{Q'}(\bm)$.
Therefore the above sum is equal to
\[
\sum_{\ov{\bm\in\Zs^{Q'_1}}{\dd_{Q'}(\bm)=0}}P'(\bm)
\prod_{a\in Q'_1}\rho(m_a)
\]
for a new polynomial $P'$.
This sum belongs to $\QM$ by \cref{QM for rho}.
\end{proof}

\begin{example}
Consider a quiver with two arrows $a,b:1\to 2$ and $f_a=f_b=\Phi_1$.
Then
\[\sum_{\ov{n_a,n_b\in\Zs}{n_a+n_b=0}}
\what f_a(n_a)\what f_b(n_b)
=\sum_{n\ne0}\frac{1}{(1-q^n)(1-q^{-n})}
=\sum_{n>0}\frac{-2q^n}{(1-q^n)^2}
=-2\sum_{n>0}Q_2(q^n)\in\QM.
\]
\end{example}

The above theorem can be interpreted as a generalization of the results in \GM.
More precisely, let $\Ga=(V,E,s)$ be a graph (see \cref{usual graphs})
with $V=[r]$ and let $\bfk\in\bN^E$.
For an orientation $Q$ of $\Ga$ (a quiver structure on $\Ga$)
define (\cf \GM[\S6])
\eq{
S(Q,\bfk)
=\sum_{\ov{w\in\Ns^E}{\dd_Q(w)=0}}
\sum_{h\in\cN(Q)}\prod_{e\in E}w_e^{k_e+1}q^{h_e w_e},
}
where $\dd_Q:\bZ^{Q_1}\to \bZ^{Q_0}$ is the boundary operator and
\begin{gather*}
\cN(Q)=\sets{h\in\bN^E}{h_a>0\ \forall a\in E_-(Q)},\\
E_-(Q)=\sets{a\in E}{s_Q(a)\ge t_Q(a)},\qquad E_+(Q)=E\ms E_-(Q).
\end{gather*}
We have $\sum_{h\ge0}q^h=\frac1{1-q}$ and
$\sum_{h>0}q^h=\frac q{1-q}$, hence
\eq{
S(Q,\bfk)
=\sum_{\ov{w\in\Ns^E}{\dd_Q(w)=0}}
\rbigg{\prod_{e\in E_+(Q)}\frac{w_e^{k_e+1}}{1-q^{w_e}}}
\rbigg{\prod_{e\in E_-(Q)}\frac{w_e^{k_e+1}q^{w_e}}{1-q^{w_e}}}
}
Define the graph sum to be the sum over all orientations of $\Ga$
\eq{S(\Ga,\bfk)=\sum_{Q\in\Or(\Ga)}S(Q,\bfk).
}
It is proved in \GM[Theorem 6.1] that if $k_e$ are even,
then $S(\Ga,\bfk)\in\QM$.
Let us interpret and prove this statement using \cref{QM acyclic}.

First of all, the contribution of every loop $e\in E$ is
\[\sum_{w,h>0}w^{k_e+1}q^{hw}
=\sum_{h>0}\psi^h\rbr{D_q^{k_e+1}\frac q{1-q}}
=\sum_{h>0}Q_{k_e+2}(q^h)=F_{k_e+2}(q),
\]
where $\psi^h(f)=f(q^h)$.
We have $F_{k_e+2}(q)\in\QM$, since $k_e+2$ is even.

Therefore we can assume that $\Ga$ has no loops.
We choose the reference orientation $Q^\circ$ of $\Ga$ such that $s(a)<t(a)$ for all $a\in Q^\circ_1$,
and let $\dd=\dd_{Q^\circ}:\bZ^E\to\bZ^V$.
The quiver $Q^\circ$ is acyclic.

\begin{lemma}
Let $\bfk\in\bN^E$ have even coordinates.
Then
\[S(\Ga,\bfk)
=\sum_{\ov{\bn\in\Zs^E}{\dd(\bn)=0}}\prod_{e\in E}
n_e^{k_e+1}\frac1{1-q^{n_e}}.
\]
In particular, $S(\Ga,\bfk)\in\QM$ by \cref{QM acyclic}.
\end{lemma}

\begin{proof}
Given $\bn\in\Zs^{E}$,
let $Q$ be the quiver obtained from $Q^\circ$ by reversing the arrows
from $E_-(\bn)=\sets{e\in E}{n_e<0}$,
and let $w=(\n{n_e})_{e\in E}\in\pN^E$.
Then $\dd(\bn)=\dd_Q(w)$ and
\[\cN(Q)=\sets{h\in\bN^E}{h_a>0\text{ if }n_a<0}.\]
This gives a bijection between $\bn\in\Zs^E$ satisfying $\dd(\bn)=0$
and the set of pairs $(Q,w)$, where $Q\in\Or(\Ga)$
and $w\in\pN^E$ satisfies $\dd_Q(w)=0$.
Let $n\in\Zs$ and $w=\n n$.
For even $k\ge0$, we have
\begin{gather*}
w^{k+1}\frac1{1-q^w}=n^{k+1}\frac1{1-q^n},\qquad n>0,\\
w^{k+1}\frac{q^w}{1-q^w}
=(-n)^{k+1}\frac{q^{-n}}{1-q^{-n}}=n^{k+1}\frac1{1-q^n},\qquad n<0.
\end{gather*}
Therefore
\[S(\Ga,\bfk)
=\sum_{Q\in\Or(\Ga)}\sum_{\ov{w\in\Ns^E}{\dd_Q(w)=0}}
\rbigg{\prod_{e\in E_+(Q)}\frac{w_e^{k_e+1}}{1-q^{w_e}}}
\rbigg{\prod_{e\in E_-(Q)}\frac{w_e^{k_e+1}q^{w_e}}{1-q^{w_e}}}
=\sum_{\ov{\bn\in\Zs^E}{\dd(\bn)=0}}\prod_{e\in E}
\frac{n_e^{k_e+1}}{1-q^{n_e}}.
\qedhere
\]
\end{proof}

\subsection{Fourier coefficients}
To study more general graph sums, we need to understand individual Fourier coefficients of quasi-elliptic functions, rather than only their constant terms.
We show that these coefficients admit explicit descriptions in terms of Lambert series, with a dependence on the Fourier indices governed by a hyperplane arrangement.

For $f\in\cC_r$ and $\bn=(n_1,\dots,n_r)\in\bZ^r$, define
\[
\what f(\bn)
=[z_1^{n_1}\dots z_r^{n_r}]\io_r(f).
\]

Let $V=\bR^r$ and let $x_i\in V^*$ for $i\in[r]$ be the corresponding coordinates and $x_I=\sum_{i\in I}x_i\in V^*$ for $I\sbs [r]$.
Let
\[\cL_r=\sets{x_I}{\es\ne I\sbs[r]}\sbs V^*\]
and let $\cA_r=\setc{\Ker(\al)\sbs V}{\al\in\cL_r}$ be the corresponding hyperplane arrangement.
A face of $\cA_r$ (or $\cL_r$) is a nonempty cone of the form
\[C_\si=\sets{v\in V}{\sgn \al(v)=\si(\al)\ \forall\al\in\cL_r}\]
for a map $\si:\cL_r\to\set{0,1,-1}$.
For a face $C$, put
\[
\cL_C=\setc{\al\in\cL_r}{\al|_C\ne0}.
\]

\begin{theorem}
\label{fourier}
Let $f\in\cC_r$ and let $C$ be a face of $\cL_r$.
Then there is a finite expression
\[
\what f(\bn)
=\sum_\nu a_\nu(q)P_\nu(\bn)
\prod_{\al\in\cL_C}f_{\nu,\al}\rbr{q^{\n{\al(\bn)}}},
\qquad
\bn\in C\cap\bZ^r,
\]
where $a_\nu\in\cC_0=\QM$, $P_\nu\in\bQ[n_1,\dots,n_r]$
and $f_{\nu,\al}\in\RR=\bQ[\frac q{1-q}]$ (\cf \cref{RR}).
\end{theorem}

\begin{proof}
We argue by induction on $r$.
For $r=0$, the statement is clear.
Let $r\ge1$ and assume the result for $\cC_{r-1}$.
Let $f\in\cC_r$ and let $C$ be a face of $\cA_r$.

Assume that $n_r=0$ on $C$.
We have
$[z_r^0]\io_r(f)=\io_{r-1}(g)$ for some $g\in\cC_{r-1}$.
Consequently, $\what f(\bm,0)=\what g(\bm)$ for $\bm\in\bZ^{r-1}$.
The restriction of $\cA_r$ to $n_r=0$ is $\cA_{r-1}$, so the assertion follows by induction on $r$.

\def\zz{u}
Assume that $n_r\ne 0$ on $C$.
For $i\in[r]$, consider the operator
\[
T_i(f)=f(z_1,\dots,qz_i,\dots,z_r)
\]
which preserves $\cC_r$.
Let $z_0=1$ and $(\zz,q)\in U_{r-1}$ with $\zz=(z_1,\dots,z_{r-1})$,
so that
\[1>\n{z_1}>\dots>\n{z_{r-1}}>\n q>0.\]
Choose $\n q<R<\n{z_{r-1}}$.
Then
\[
f_n(\zz):=[z_r^n] \io_r(f)
=\frac1{2\pi i}\int_{\n x=R}\om(x),\qquad
\om(x)=f(\zz,x)x^{-n-1}dx,
\]
and similarly
\[
[T_rf]_n(\zz)
=\frac1{2\pi i}\int_{\n x=R}f(\zz,qx)x^{-n-1}dx
=\frac{q^n}{2\pi i}\int_{\n x=\n qR}\om(x).
\]

The poles of $f(\zz,x)$ in the annulus $\n qR\le \n x\le R$
are contained in $\sets{qz_i}{0\le i<r}$ (\cf~\cite{hongminh_quasia}).
By the residue theorem, we obtain
\[
g_n(\zz):=[q^nf-T_rf]_n(\zz)
=q^n\sum_{i=0}^{r-1}\Res_{x=qz_i}\om(x).
\]
Using $h=(1-T_r)f$, we obtain
\[f_n=\frac1{1-q^n}(h_n-g_n),\qquad n\ne0.\]
Let $\bn=(\bm,n)\in C\cap\bZ^r$, where $\bm=(n_1,\dots,n_{r-1})$ and $n=n_r\ne0$.
Then
\eq{\label{fourier:recursion}
\what f(\bn)=[\zz^{\bm}]f_n
=\frac1{1-q^n}\rbr{\what h(\bn)-[u^\bm]g_n}
}

Let us analyze $[u^\bm]g_n$.
At the possible pole $x=qz_i$ of $f(u,x)$,
write $x=qz_i e^s$.
Then $dx/x=ds$ and
$q^nx^{-n}=z_i^{-n}e^{-ns}$.
By~\cref{Phi:local},
we have $f(u,qz_ie^s)=(T_rf)(u,z_ie^s)\in\cC_{r-1}\lser s$.
Therefore
\[
f(\zz,qz_ie^s)
=\sum_{\ell\ge-M_i}c_{i,\ell}(\zz)s^\ell,
\]
where $c_{i,\ell}\in\cC_{r-1}$.
Then
\[q^n\Res_{x=qz_i}
f(\zz,x)x^{-n}\frac{dx}{x}
=z_i^{-n}\Res_{s=0}f(\zz,qz_ie^s)e^{-ns}\,ds
=z_i^{-n}\sum_{k=0}^{M_i-1}\frac{(-n)^k}{k!}c_{i,-k-1}(\zz).
\]
Consequently,
\[
g_n(\zz)=\sum_{i=0}^{r-1}z_i^{-n}
\sum_\nu P_{i,\nu}(n)g_{i,\nu}(\zz),
\]
where
$P_{i,\nu}\in\bQ[n]$ and $g_{i,\nu}\in\cC_{r-1}$.
Let $e_1,\dots,e_{r-1}$ be the standard basis of $\bZ^{r-1}$ and put $e_0=0$.
Then
\begin{equation}
\label{fourier:residue}
[\zz^{\bm}]g_n(\zz)
=
\sum_{i=0}^{r-1}\sum_\nu
P_{i,\nu}(n)
\what g_{i,\nu}(\bm+ne_i).
\end{equation}

For $0\le i\le r-1$,
consider the linear map
\[\pi_i:\bR^r\to\bR^{r-1},\qquad
(\bm,n)\mto \bm+n e_i.
\]
For every $\es\ne J\sbs[r-1]$, we have
\[
x_J\circ\pi_i
=\begin{cases}
x_J&i\notin J,\\
x_{J\cup\set r}&i\in J.
\end{cases}
\]
Therefore $x_J$ has constant sign on $\pi_i(C)$ and $\pi_i(C)$
is contained in a unique face $C_i$ of $\cA_{r-1}$.
By induction on $r$, the Fourier coefficients
$\what g_{i,\nu}$ admit the required expansion on $C_i$.
For every linear form $x_J\in\cL_{C_i}$ occurring in this expansion,
the form $x_I=x_J\circ\pi_i$
has constant nonzero sign on~$C$, hence $x_I\in\cL_C$.
Therefore the right-hand side of \cref{fourier:residue}
has the required expansion on $C$.

Now let us analyze $\what h(\bn)$.
Following \cite{hongminh_quasia},
let $B\sbs\cC_r$ be the subalgebra generated by $\cC_0$
and
\begin{enumerate}
\item $\Phi_k(z_j/z_i)$ for $0\le i\ne j\le r$ and $k\ge2$,
\item $\Phi_1(z_r/z_i)-\Phi_1(z_r)$ for $i\ne r$,
\item $\Phi_1(z_j/z_i)$ for $0\le i\ne j<r$.
\end{enumerate}
The elements of $B$ are $T_r$-invariant and
$\cC_r=B[\ga]$, where $\ga=\Phi_1(z_r)$
satisfies $T_r \ga=\ga+1$.

We now use induction on the degree of $f$ in $\ga$.
For $h=(1-T_r)f$, we have $h=0$ if $\deg_\ga f=0$,
and $\deg_\ga h<\deg_\ga f$ otherwise.
Thus, by induction, $\what h(\bn)$ has the required expansion on $C$.

Since $n_r\ne0$ on $C$, we have $x_{r}\in\cL_C$.
Moreover,
\[
\frac1{1-q^n}
=
\begin{cases}
1+Q_1(q^{|n|}),&n>0,\\
-Q_1(q^{|n|}),&n<0,
\end{cases}
\]
where $Q_1(q)=\frac q{1-q}\in\RR$.
Since $\RR$ is an algebra, \cref{fourier:recursion}
implies that $\what f(\bn)$ has the required expansion on $C$.
This completes the induction.
\end{proof}

\section{Multiple \tpdf{$q$}{q}-zeta values and bi-brackets}
\label{mqzv:sec}
We now turn to the algebra of multiple $q$-zeta values in which our generating series will lie.
Its description in terms of bi-brackets provides a natural framework for summing the Fourier coefficients obtained in the preceding section.
We establish general summation results that will allow us to recognize certain graph sums as elements of this algebra.

\subsection{Multiple \tpdf{$q$}{q}-zeta values}
Let $A=\bQ\pser q$ be equipped with the \la-ring structure with Adams operations $\psi^n(f)=f(q^n)$ for $f\in A$ and $n\ge1$
(see \eg \cite{atiyah_group}).
For $f_1,\dots,f_r\in A$ with $f_1\in A^+=qA$,
define the multiple $q$-zeta value
\eq{Z(f_1\ts\dots \ts f_r)
=\sum_{n_1>\dots>n_r>0}\prod_{i=1}^r\psi^{n_i}(f_i)\in A.
}
This defines an algebra morphism
\eq{
Z:\QS(A^+)\to A,
}
where $\QS(A^+)=\bop_{r\ge0}(A^+)^{\ts r}$
is the quasi-shuffle algebra of the non-unital algebra $A^+$
(see \eg \cite{mozgovoy_multiple}).
For a subalgebra $B\sbs A$,
let $B^+=B\cap A^+$ and consider the subalgebra
\eq{\QS(B^+,B)=\bQ\oplus\rbr{B^+\ts\QS(B)}\sbs\QS(B).
}
Then $Z$ extends to an algebra morphism $Z:\QS(B^+,B)\to A$.

Recall that $Q_k(q)=D_q^{k-1}\frac q{1-q}$ for $k\ge1$ and $Q_0(q)=1$.
Consider the algebra (\cf \cite{mozgovoy_multiple})
\eq{\label{RR}
\RR
=\sets{g/(1-q)^k}{g\in\bQ[q],\, \deg g\le k}
=\bQ[Q_1]=\bop_{k\ge0}\bQ Q_k.}
Then $\RR^+=\RR\cap q\bQ\pser q=\bop_{k\ge1}\bQ Q_k$.
Define the algebras of multiple $q$-zeta values
\eq{
\Zq^+=\Im\rbr{Z:\QS(\RR^+)\to\bQ\pser q},\qquad
\Zq=\Im\rbr{Z:\QS(\RR^+,\RR)\to\bQ\pser q}.
}
It was conjectured in \cite{bachmann_dimension}
and proved in \cite{hirose_unified} that
\eq{\Zq^+=\Zq.}
Since $F_k(q)=\sum_{d\ge1}Q_k(q^d)=Z(Q_k)\in\cZ_q$ for $k\ge1$,
we conclude that
\[\QM\sbs\Zq.\]

\begin{remark}
Let $f_{1,k}=\frac{q^{k-1}}{(1-q)^k}$ for $k\ge1$ and let
(\cf \cite{mozgovoy_multiple})
\[
\RR_1=\sets{g/(1-q)^k}{g\in\bQ[q],\, \deg g<k}
=\bop_{k\ge1}\bQ f_{1,k}.\]
The algebra
$\cZ_{q,1}=\Im\rbr{Z:\QS(\RR_1^+)\to\bQ\pser q}$
is also sometimes called the algebra of multiple $q$-zeta values
(see \eg \cite{okounkov_hilbert}, where it is denoted
by $\opn{qMZV}$).
Strictly speaking, the conjecture of Okounkov \cite{okounkov_hilbert}
concerns the algebra~$\cZ_{q,1}$,
whereas we work with the algebra $\Zq$.
\end{remark}

\subsection{Bi-brackets}
For $s\in\pN^r$ and $t\in\bN^r$, define the bi-bracket
(\cf \cite{bachmann_algebraa})
\eq{\label{bb}
\bb st
=\sum_{n_1>\dots>n_r>0}\sum_{d\in\pN^r}
\prod_{i=1}^r n_i^{t_i}d_i^{s_i-1}q^{d_in_i}
=\sum_{n_1>\dots>n_r>0}
\prod_{i=1}^r n_i^{t_i}Q_{s_i}(q^{n_i})\in\bQ\pser q.
}
Let $\BD\sbs \bQ\pser q$ be the subspace spanned by all bi-brackets (and $1$).
It is proved in \cite[Theorem 3.6]{bachmann_algebraa}
that~$\BD$ is a subalgebra of $\bQ\pser q$.
We give an alternative argument.

\begin{lemma}\label{BD alg}
The subspace $\BD\sbs\bQ\pser q$ is a subalgebra.
\end{lemma}
\begin{proof}
We equip the algebra $\bQ\pser q[y]$ with the \la-ring structure
with Adams operations
\[\psi^n(f(q)y^k)=f(q^n)n^k y^k.
\]
For $\BB^+=\RR^+[y]$, there is an algebra morphism
\eq{\label{Z BDy}
Z:\QS(\BB^+)\to\bQ\pser q[y],\qquad
Z(f_1\ts\dots\ts f_r)=\sum_{n_1>\dots>n_r>0}
\prod_{i=1}^r\psi^{n_i}(f_i).
}
In particular, for $f_i=Q_{s_i}(q)y^{t_i}=\sum_{d\ge1}d^{s_i-1}q^dy^{t_i}$,
we obtain
\eq{
Z(f_1\ts\dots\ts f_r)
=\sum_{n_1>\dots>n_r>0}\sum_{d\in\pN^r}
\prod_{i=1}^rd_i^{s_i-1}q^{d_in_i}n_i^{t_i}y^{t_i}
=y^{\sum_it_i}\bb st.
}
The elements $Q_sy^t$ for $s\ge1$ and $t\ge0$ form a basis of $\BB^+$.
Therefore the image $\BD_y\sbs\bQ\pser q[y]$
of the algebra morphism \cref{Z BDy}
is spanned by the elements $y^{\sum_it_i}\bb st$.
Its specialization $\BD\sbs\bQ\pser q$ at $y=1$
is a subalgebra spanned by the elements $\bb st$.
\end{proof}

The following result is proved in
\cite[Theorem 2.3]{bachmann_dimension}.
We give an alternative argument.

\begin{lemma}\label{Z=BD}
We have $\Zq=\BD$.
\end{lemma}

\begin{proof}
The inclusion $\Zq\sbs\BD$ follows from \cref{BD1:lm}.
Conversely, consider a bi-bracket~\cref{bb}.
The polynomials
\[
p_m(n)=\binom{n-1}{m}
=\sum_{n>k_1>\dots>k_m>0}1
\]
for $m\ge0$, form a basis of $\bQ[n]$.
Hence, by linearity,
it is enough to replace every $n_i^{t_i}$ in a bi-bracket $\bb st$
by $p_{m_i}(n_i)$.
We therefore have to consider (for $m_1,\dots,m_r\ge0$)
\[
\sum_{n_{1}>\dots>n_{r}>0}
\prod_{i=1}^r
\rbigg{Q_{s_i}(q^{n_i})
\sum_{n_{i}>n_{i,1}>\dots>n_{i,m_i}>0}1}.
\]
Partition the summation domain according to the weak order of
the variables $n_i$ and $n_{i,j}$.
On each stratum, let
$N_1>\dots>N_k>0$ be their distinct values.
The contribution of the stratum is
\[
Z(f_1,\dots,f_k)
=
\sum_{N_1>\dots>N_k>0}\prod_{a=1}^k f_a(q^{N_a}),
\]
where every $f_a$ is either $1$ or one of the $Q_{s_i}$.
Moreover, $N_1=n_1$ and $f_1=Q_{s_1}\in\RR^+$.
Hence the contribution of every stratum belongs to $\Zq$.
Therefore $\BD\sbs\Zq$.
\end{proof}

\subsection{Closure properties of bi-brackets}
We establish several closure properties of the algebra of bi-brackets that will be useful for evaluating graph sums.
In particular, we show that polynomial weights and certain partial ordering constraints on summation variables preserve membership in this algebra.

\begin{lemma}\label{Faulhaber:lm}
Let $P\in\bQ[n_1,\dots,n_r]$.
Then
\[\sum_{n\ge n_1>\dots>n_r>m}P(n_1,\dots,n_r)\]
is a polynomial in $n\ge m\ge 0$.
\end{lemma}
\begin{proof}
We can assume that
$P(\bn)=\prod_{i=1}^r n_i^{t_i}$ for some
$t=(t_1,\dots,t_r)\in\bN^r$.
Let $f_t(m,n)=\sum_{n\ge n_1>\dots>n_r>m}\prod_{i=1}^r n_i^{t_i}$.
By Faulhaber's formula, for $p\ge0$,
\[
g_p(n)=\sum_{i=1}^n i^p=\frac{1}{p+1}\sum_{k=0}^p\binom{p+1}k B^+_k n^{p+1-k},
\]
where $B^+_k$ are Bernoulli numbers defined by
$\sum_{k\ge0}B^+_kt^k/k!=\frac{te^t}{e^{t}-1}$.
In particular, $g_p(n)$ is a polynomial in $n\ge 0$.
For $r=1$, we have $f_t(m,n)=g_t(n)-g_t(m)$,
which is a polynomial in~$m,n$.
By induction on $r$,
\[
f_{t_1,\dots,t_{r+1}}(m,n)
=\sum_{j=m+1}^{n}j^{t_{r+1}} f_t(j,n),
\]
where $t=(t_1,\dots,t_r)$ and $f_t(j,n)$ is a polynomial in $j,n$.
By Faulhaber's formula, the sum is a polynomial in $m,n$.
\end{proof}

\begin{theorem}\label{BD1:lm}
Let $P\in\bQ[n_1,\dots,n_r]$ and $f_1,\dots,f_r\in \RR$ with $f_1\in \RR^+$.
Then
\[
\sum_{n_1>\dots>n_r>0}P(\bn)\prod_{i=1}^r f_i(q^{n_i})\in\BD.
\]
\end{theorem}
\begin{proof}
We can assume that $P(\bn)=\prod_{i=1}^r n_i^{t_i}$ and $f_i=Q_{s_i}$ for some $t_i,s_i\ge0$ with $s_1>0$.
If all $s_i>0$,
then the sum is equal to
\[
\sum_{n_1>\dots>n_r>0}\prod_{i=1}^r n_i^{t_i}Q_{s_i}(q^{n_i})
=\bb st\in\BD.
\]
Generally, let
\[
J=\sets{i\in[r]}{s_i>0}=\set{j_1<\dots<j_k}.
\]
Since $s_1>0$, we have $j_1=1$.
If $n_1>\dots>n_r>0$, then the numbers $N_a=n_{j_a}$ satisfy
\[
N_1>\dots>N_k> N_{k+1}=0.
\]

For fixed $N_1>\dots>N_k>0$,
sum first over all $n_i$ with $i\notin J$.
The variables $n_i$ with $i\notin J$ occur only in the polynomial factor, since $Q_{s_i}=Q_0=1$.
Therefore the sum is a polynomial $\wtl P\in\bQ[N_1,\dots,N_k]$
by \cref{Faulhaber:lm}.
Hence the original sum is equal to
\[
\sum_{N_1>\dots>N_k>0}
\wtl P(N_1,\dots,N_k)
\prod_{a=1}^k Q_{s_{j_a}}(q^{N_a}),
\]
with $s_{j_a}>0$.
It belongs to $\BD$ by the previous argument.
\end{proof}

Given posets $S$ and $T$, let
\[\cI(S,T)=\sets{\al:S\to T}{\al(i)<\al(j)\text{ for }i<j}\sbs T^S
\]
be the set of strictly increasing maps.

\begin{theorem}\label{BD2:th}
Let $S$ be a finite poset, $P(\bn)\in\bQ[n_i\col i\in S]$,
and $f_i\in\RR$ for $i\in S$ such that $f_i\in\RR^+$
for all maximal $i\in S$.
Then
\[
\sum_{\bn\in\cI(S,\pN)}
P(\bn)\prod_{i\in S}f_i(q^{n_i})
\in\BD.
\]
\end{theorem}

\begin{proof}
We can assume that $P(\bn)=\prod_{i\in S}n_i^{t_i}$ for some $t_i\ge0$.
For every $\bn\in\cI(S,\pN)$,
there exists a unique ordered set partition
$S=B_1\sqcup\dots\sqcup B_r$
and integers $N_1>\dots>N_r>0$
such that $n_i=N_k$ for $i\in B_k$.
For $B=(B_1,\dots, B_r)$,
the stratum $\cI_B\sbs \cI(S,\pN)$ is non-empty \ifft
$i\in B_k$ and $j\in B_l$ with $i<j$ implies $k>l$.
The contribution of the stratum is
\[
\sum_{N_1>\cdots>N_r>0}
P_B(N_1,\ldots,N_r)\prod_{k=1}^r g_k(q^{N_k}),
\]
where $P_B=\prod_{k=1}^r\prod_{i\in B_k}N_k^{t_i} \in\bQ[N_1,\ldots,N_r]$
and $g_k=\prod_{i\in B_k}f_i\in\RR$.
If $i\in B_1$, $j\in B_l$ and $i<j$, then $1>l$, a contradiction.
Therefore all elements of $B_1$ are maximal and $g_1\in\RR^+$.
Therefore the contribution belongs to $\BD$ by \cref{BD1:lm}.
\end{proof}

\begin{conjecture}
\label{conj1}
Let $Q$ be a quiver with a decomposition
$Q_1=E_+\sqcup E_-$ such that reversing the arrows from $E_-$ gives an acyclic quiver.
Let $P\in\bQ[n_a\col a\in Q_1]$.
Then
\[\sum_{\ov{\bn\in\Ns^{Q_1}}{\dd(\bn)=0}}
P(\bn)\prod_{a\in E_+}\rho(n_a)\prod_{a\in E_-}\rho(-n_a)\in\BD,
\]
where $\rho(n)=\frac1{1-q^n}=\sum_{k\ge0}q^{kn}$
and $\rho(-n)=\frac1{1-q^{-n}}=-\sum_{k\ge1}q^{kn}$
for $n>0$.
\end{conjecture}

This conjecture should be compared to \cref{QM for E_-}.
A closely related statement appears in
\cite[Prop.~5.3]{qin_multiple}, but the proof of
\cite[Lemma~5.2]{qin_multiple}, on which it relies,
is incomplete.
The following is a stronger conjecture.

\begin{conjecture}
Let $Q$ be a quiver and let
$P\in\bQ[n_a\col a\in Q_1]$
and $f_a\in\RR$ for $a\in Q_1$ be such that
$Q$ restricted to $\sets{a\in Q_1}{f_a\notin\RR^+}$ is acyclic.
Then
\[\sum_{\ov{\bn\in\Ns^{Q_1}}{\dd(\bn)=0}}
P(\bn)\prod_{a\in Q_1}f_a(q^{n_a})\in\BD,
\]
\end{conjecture}

\subsection{Nested families}
The Fourier expansions arising in our applications lead to Lambert series whose arguments involve sums over subsets of variables.
We show that when these subsets form a nested family, the resulting sums belong to the algebra of bi-brackets, provided a suitable finiteness condition holds.

\begin{definition}
A family $\cS\sbs 2^{[r]}$ of nonempty subsets of $[r]$
is called nested (or laminar) if for all $I,J\in\cS$
we have $I\sbs J$, $J\sbs I$ or $I\cap J=\es$.
\end{definition}

We order a nested family $\cS\sbs 2^{[r]}$ by inclusion.
Then
\[\cS_{\ge A}=\sets{I\in\cS}{I\sps A}\]
is a chain for every nonempty subset $A\sbs [r]$.
Therefore the Hasse diagram of $\cS$ is a forest,
with the maximal elements of $\cS$ being the roots of the trees.

\begin{remark}
It is common to require a nested family $\cS\sbs 2^{[r]}$
to contain only $I\in\cS$ with $\n I\ge2$
(see \eg \cite{fulton_compactification}).
Such a family can be identified with a family that contains all singletons.
\end{remark}

\begin{theorem}
\label{lem:proper-laminar}
Let $P\in\bQ[n_1,\ldots,n_r]$
and let $\cS\sbs2^{[r]}$ be a nested family
and $f_I\in\RR$ for $I\in\cS$ be such that
$\cS^+=\sets{I\in\cS}{f_I\in\RR^+}$
satisfies $\bigcup\cS^+=\bigcup_{I\in\cS^+}I=[r]$.
Then
\[
\sum_{\bn\in\pN^r}P(\bn)\prod_{I\in\cS}f_I(q^{n_I})\in\BD,
\]
where $n_I=\sum_{i\in I}n_i$ for $I\sbs[r]$.
\end{theorem}

\begin{proof}
Let $f_I=c_I+f_I^+$, where $c_I\in \bQ$ and $f_I^+\in\RR^+$.
Then
\[
\sum_{\bn\in\pN^r}P(\bn)\prod_{I\in\cS}f_I(q^{n_I})
=\sum_{\cS^+\sbs\cT\sbs\cS}
\rbr{\prods_{I\in\cS\ms\cT}c_I}
\sum_{\bn\in\pN^r}P(\bn)\prod_{I\in\cT}f_I^+(q^{n_I}).
\]
Every $\cT$ above is nested and satisfies $\bigcup\cT=\bigcup\cS^+=[r]$.
Therefore we can assume that $f_I\in\RR^+$ for all $I\in\cS$.
We can assume that $P(\bn)=\prod_{i=1}^r n_i^{t_i}$
and $f_I=Q_{s_I}$ for some $t_i\ge0$ and $s_I\ge1$.
Let $J_1,\ldots,J_c$ be the maximal elements of $\cS$.
Then $[r]=J_1\sqcup\dots\sqcup J_c$.
Let $\cS_k=\sets{I\in\cS}{I\sbs J_k}$ for $k\in[c]$.
Then the required sum factors as
\[
\Phi=\prod_{k=1}^c
\rbr{\sum_{\bn\in\pN^{J_k}}
\prod_{i\in J_k}n_i^{t_i}\prod_{I\in\cS_k}Q_{s_I}(q^{n_I})
}.
\]
Therefore we can assume that $\cS$ has a unique maximal member $J_1=[r]$,
and the Hasse diagram of $\cS$ is a rooted tree with the root $[r]$.
Then
\[
\Phi=
\sum_{\bn\in\pN^r}\prod_{i\in[r]}n_i^{t_i}
\prod_{I\in\cS}Q_{s_I}(q^{n_I})
=\sum_{\substack{\bn\in\pN^r\\ d_I>0,\, I\in\cS}}
\rbigg{\prod_{i\in[r]} n_i^{t_i}}
\rbigg{\prod_{I\in\cS}d_I^{s_I-1}}
q^{\sum_{I\in\cS}d_In_I}.
\]
Given $d\in\pN^\cS$,
define $N_I=\sum_{J\in\cS_{\ge I}}d_J$ for $I\in\cS$.
Then the change of variables $d\mapsto N$ is a bijection between $\pN^\cS$ and the set $\cI(\cS^\op,\pN)$
of tuples $N\in\pN^\cS$
satisfying $N_I>N_J$ for $I<J$ (meaning that $I\subsetneq J$).

For each $i\in[r]$, the set $\cS_i=\sets{I\in\cS}{I\ni i}$
is a nonempty chain.
Let $I_i\in\cS_i$ be its minimal member.
Then $\sum_{I\in\cS_i}d_I=N_{I_i}$, hence
\[
\sum_{I\in\cS}d_In_I
=\sum_{I\in\cS}\sum_{i\in I}d_In_i
=\sum_{i\in[r]} n_iN_{I_i}.
\]
This implies
\[\Phi
=\sum_{N\in\cI(\cS^\op,\pN)}
\rbigg{\prod_{I\in\cS} d_I^{s_I-1}}
\prod_{i\in[r]}\rbr{\sum_{n>0}n^{t_i}q^{nN_{I_i}}}
=\sum_{N\in\cI(\cS^\op,\pN)}
\rbigg{\prod_{I\in\cS} d_I^{s_I-1}}
\prod_{i\in[r]}Q_{t_i+1}\rbr{q^{N_{I_i}}}
\]
Therefore
\[
\Phi
=\sum_{N\in\cI(\cS^\op,\pN)}
\wtl P(N)
\prod_{I\in\cS}g_I(q^{N_I}),
\]
where $\wtl P$ is a polynomial and
$g_I(q)=\prod_{i:I_i=I}Q_{t_i+1}(q)\in\RR$ for $I\in\cS$.
If $I\in\cS$ is a minimal element, then $I_i=I$ for all $i\in I$, hence $g_I\in\RR^+$.
By \cref{BD2:th} applied to the opposite poset~$\cS^\op$,
we obtain $\Phi\in\BD$.
\end{proof}

We say that a map $F:X\to\bQ\pser q$ is $q$-proper if,
for every $d\ge0$, the set
\[\supp([q^d]F)=\sets{x\in X}{[q^d]F(x)\ne0}\]
is finite.
Then $\sum_{x\in X}F(x)\in\bQ\pser q$ is well-defined.

\begin{theorem}[Proper summation]
\label{proper sum}
Let $F:\pN^r\to\bQ\pser q$ be a $q$-proper map
such that $F(\bn)$ is a finite $\QM$-linear combination of terms
\[
P(\bn)\prod_{I\in\cS}f_I(q^{n_I}),
\]
where $P\in\bQ[n_1,\dots,n_r]$, $\cS\sbs 2^{[r]}$ is a nested family and $f_I\in\RR$ for $I\in\cS$.
Then
\[
\sum_{\bn\in\pN^r}F(\bn)\in\BD.
\]
\end{theorem}

\begin{proof}
Write every $f_I$ as $f_I(0)+f_I^+$ with $f_I^+\in\RR^+$.
Expanding $F$, we conclude that $F(\bn)$
is a finite sum of terms
\[T_k(\bn)=a_k(q)P_k(\bn)\prods_{I\in\cS_k}f_{k,I}(q^{n_I}),\]
where $a_k(q)\in\QM$, $P_k\in\bQ[n_1,\dots,n_r]$,
$\cS_k$ is a nested family and $f_{k,I}\in\RR^+$.

Let us show that the contribution of terms for which
$\bigcup\cS_k\ne[r]$ vanishes.
For a proper subset $S\sbs[r]$,
let $F_S=\sum_{\bigcup\cS_k=S}T_k$ and
$F_{\le S}=\sum_{\bigcup\cS_k\sbs S}T_k$.
Let $\bn\in\pN^r$ and $d\ge0$ be fixed and let
$\bv=\one_{[r]\ms S}\in\bN^r$.

If $\bigcup\cS_k\not\sbs S$,
then $I\not\sbs S$ for some $I\in\cS_k$.
Therefore $n_I+tv_I\to\infty$ for $t\to\infty$.
Since $f_{k,I}\in\RR^+$, we obtain
$[q^d]T_k(\bn+t\bv)=0$ for $t\gg0$.
Therefore
\[
[q^d]F(\bn+t\bv)=[q^d]F_{\le S}(\bn+t\bv)
\]
for $t\gg0$.
The left-hand side is eventually zero by $q$-properness, while the right-hand side is a polynomial in $t$.
Therefore it vanishes identically, hence $[q^d]F_{\le S}(\bn)=0$.
This implies that $F_{\le S}=0$.
Then $F_S=0$ by induction on $\n S$, hence every nonproper contribution vanishes.

This implies that we should consider only the terms with $\bigcup\cS_k=[r]$.
Such terms are contained in $\BD$ by \cref{lem:proper-laminar}
and the fact that $a_k(q)\in \QM\sbs\BD$.
\end{proof}

\section{Wick formulas}
\label{wick:sec}
Wick's formalism provides the link between the operator-theoretic and combinatorial aspects of our problem.
We develop formulas expressing normalized traces of Heisenberg operators as sums over graphs, with contractions governed by the underlying Frobenius algebra.
These formulas will allow us to exploit the special properties of the cohomology of a surface.

\subsection{Graphs and matchings}
A set partition of a finite set $I$ is a collection $\pi\sbs 2^I$
of nonempty subsets such that $I=\bigsqcup_{B\in\pi}B$.
It can be identified with an equivalence relation on $I$.
Let $\pi^{(r)}=\sets{B\in\pi}{\n B=r}$ for $r\ge1$.
Let $\Pi_I$ denote the set of all set partitions of $I$ and let
\[
\Pi^{\le r}_I=\sets{\pi\in\Pi_I}{\n B\le r\ \forall B\in\pi},\qquad
\Pi^{(r)}_I=\sets{\pi\in\Pi_I}{\n B=r\ \forall B\in\pi}.
\]
Define the sets of matchings (or partial pairings)
and perfect matchings (or pairings) on $I$ to be
\[\M_I=\Pi^{\le2}_I,\qquad \PM_I=\Pi^{(2)}_I.
\]
For $I=[n]$, we write $\M_n=\M_I$ and $\PM_n=\PM_I$.

\begin{definition}
\label{graph w legs}
A \idef{graph with legs}
(\cf \cite{serre_trees,getzler_modular,getzler_operadsb})
is a tuple $G=(V,F,p,\si)$,
where
\begin{enumerate}
\item $V$ and $F$ are sets, called the sets of vertices and flags,
respectively.
\item $p:F\to V$ is a map,
where $p\inv(v)$ is called the set of flags incident with $v\in V$.
\item $\si:F\to F$ is an involution,
whose fixed points are called legs and the remaining orbits are called edges.
\end{enumerate}
The involution $\si$ can be identified with the matching $\pi\in\M_F$ consisting of $\si$-orbits,
thus we also denote a graph by $(V,F,p,\pi)$.
Let $L(G)=F^\si\iso\pis$ be the set of legs and
\[E(G)=\sets{\set{f,\bar f}}{f\in F\ms F^\si}=\pip,\qquad
\bar f=\si f,\]
be the set of edges.
We usually call a graph with legs simply a graph
and a graph with $L(G)=\es$ a proper graph.
We say that an edge $e=\set{f,\bar f}\in E(G)$
connects the vertices $p(f)$ and $p(\bar f)$,
and we call $e$ a loop if $p(f)=p(\bar f)$.
\end{definition}

\begin{remark}
\label{usual graphs}
Alternatively, a graph is defined to be a triple $(V,E,s)$, where $s:E\to 2^V$ is a map such that $\n{s(e)}\in\set{1,2}$
for all $e\in E$.
A proper graph $G=(V,F,p,\si)$ can be identified with
the graph $(V,E,s)$,
where $E=E(G)$ and $s:E\to 2^V$, $e=\set{f,\bar f}\mto
\set{p (f),p(\bar f)}$.

\end{remark}

\begin{remark}
The realization $\n G$ of a graph with legs $G=(V,F,p,\si)$
is the quotient of $V\sqcup(F\xx[0,1])$
by the equivalence relation generated by
\[(f,t)=(\bar f,1-t),\qquad(f,0)=p(f),\qquad
f\in F,\,t\in[0,1].\]
We may identify $\n G$ with the quotient of $V\sqcup(F\xx[0,\oh])$
by the equivalence relation generated by
$(f,\oh)=(\bar f,\oh)$ and $(f,0)=p(f)$ for $f\in F$.
We equip $\n G$ with the quotient topology.
\end{remark}

Let $b_i(G)$ denote the $i$-th Betti number of $G$.
Thus $e(G)=b_0(G)-b_1(G)=\n{V(G)}-\n{E(G)}$.
If $G$ is connected, then $b_1(G)=\n{E(G)}-\n{V(G)}+1$.
Let $\pi_0(G)$ denote the set of connected components of $G$.
We treat every $C\in\pi_0(G)$ as a subgraph of $G$,
which can be identified with the set of vertices $V(C)\sbs V(G)$.

\subsection{Wick formulas}
Let $H$ be a quadratic space over a field $\bk$ and
let $V=V(H)$ be the Fock space of $H$.
Let $\d$ be the degree operator on $V$,
thus $\d(x)=nx$ for $x\in V_n$.
Let $A=\iend(V)$ be the graded algebra of endomorphisms.
For $f\in A$, define the $q$-trace
\[\tr_q(f)=\tr(q^\d f)=\sum_{n\ge0}q^n \tr(f_{nn}),
\]
where $f_{nn}\in\Hom(V_n,V_n)$ is the component of $f$.
Consider the normalized $q$-trace
(see \eg \cite{hongminh_quasia} for more details)
\eq{
\ta:A\to \bk\lser q,\qquad
\ta(f)=\frac{\tr(q^\d f)}{\tr(q^\d)}.}
Then $\ta(1)=1$ and for all $a\in A_n$ and $b\in A$ we have
\eq{\label{q-tr}
\ta(ab)=q^n\cdot \ta(ba).}
In particular, if $a\in A_0$, then $\ta(ab)=\ta(ba)$.

\begin{remark}[Normalized trace]
Generally, a normalized $q$-trace on a graded algebra $A$
is a linear map $\ta:A\to B$ to a commutative algebra $B$
such that $\ta(1)=1$ and \cref{q-tr} holds for a fixed invertible element $q\in B$.
\end{remark}

\begin{lemma}
Let $a\in A_n$ and $b\in A$ for $n\ne0$.
Then
\[\ta(ab)=\ta([a,b])/(1-q^{-n}).\]
\end{lemma}
\begin{proof}
We have $\ta(ab)=q^n\ta(ba)$, hence $\ta([a,b])=(1-q^{-n})\ta(ab)$.
\end{proof}

In what follows, we consider (\cf \cref{hPhi})
\eq{
\ra_k(n)
=\hPhi_k(n)=\frac{n^{k-1}}{1-q^n}\qquad k\ge1,\, n\ne0,
}
and $\ra(n)=\ra_1(n)=\frac1{1-q^n}$.
For $m,n\ne0$ and $a,b\in H$, we have $a_m\in A_{-m}$ and
\eq{\label{ta ab}
\begin{gathered}
\ta(a_mb_n)
=\de_{m,-n}(a,b)\frac{m}{1-q^m}
=\de_{m,-n}(a,b)\rb(m),\\
\ta(\nop{a_mb_n})=\de_{m,-n}(a,b)\rb(-\n m),
\end{gathered}
}
since $\ta(\nop{a_mb_n})=\ta(\nop{b_na_m})=\ta(a_mb_n)$ if $m<n$.

Consider the space of Heisenberg operators
\eq{\Ho=\Im\rbr{H[t^{\pm1}]\to\End(V)}.}

\begin{lemma}[Wick formula]
\label{q-Wick:lm}
Let $a_1,\dots,a_r\in\Ho$.
Then
\begin{equation}\label{eq:q-Wick}
\ta(a_1\dots a_r)
=\sum_{\pi\in\PM_r}\prod_{\set{i<j}\in\pi}\ta(a_ia_j).
\end{equation}

\end{lemma}
\begin{proof}
We can assume that all $a_i$ are homogeneous,
say $a_i\in\Ho_{k_i}$ with $k_i\ne0$.
The commutators $c_{ij}=[a_i,a_j]$ are scalars.
Therefore
\[[a_1,a_2\dots a_r]
=\sum_{j=2}^r c_{1j}a_2\dots \what{a_j}\dots a_r.\]
We have $\ta([a_1,b])=(1-q^{-k_1})\ta(a_1b)$, hence
\[(1-q^{-k_1})\ta(a_1\dots a_r)
=\sum_{j=2}^r c_{1j}\ta(a_2\dots \what{a_j}\dots a_r).\]
Since $\ta(a_1a_j)=c_{1j}/(1-q^{-k_1})$, we obtain
\[\ta(a_1\dots a_r)
=\sum_{j=2}^r \ta(a_1a_j)\ta(a_2\dots \what{a_j}\dots a_r).\]
Now we apply induction on $r$.
\end{proof}

\begin{lemma}\label{wick-exp1}
Let $a_1,\dots,a_r,b\in\Ho$. Then
\[
\ta(a_1\dots a_r e^{tb})
=\ta(e^{tb})\sum_{\pi\in\M_r}
\prod_{\set{i<j}\in\pi}\ta(a_i a_j)
\prod_{\set i\in\pi} t\ta(a_i b).
\]
\end{lemma}

\begin{proof}
We expand
$e^{tb}=\sum_{k\ge0}\frac{t^k}{k!}b^k$
and apply \cref{q-Wick:lm} to $a_1\dots a_r b^k$.
Fix $\pi\in\M_r$ and put $\ell=\n{\pi^{(1)}}$.
In a pairing extending $\pi$, every $\set i\in \pi$ must be paired with one of the $k$ copies of~$b$.
There are $\frac{k!}{(k-\ell)!}$
ways to choose and assign these copies.
The remaining $k-\ell$ copies of $b$ pair among themselves.
Therefore the contribution of $\pi$ is
\[\prod_{\set{i<j}\in\pi}\ta(a_i a_j)
\prod_{\set i\in \pi}\ta(a_i b)
\sum_{k\ge \ell}\frac{t^k}{k!}\frac{k!}{(k-\ell)!}\ta(b^{k-\ell})
=\ta(e^{tb})
\prod_{\set{i<j}\in\pi}\ta(a_i a_j)
\prod_{\set i\in\pi}t\ta(a_i b)
.
\]
Summing over all matchings $\pi$ proves the formula.
\end{proof}

Let $b_1,\dots,b_s\in\Ho$ and let $t_1,\dots,t_s$ be variables.
If $c=\oh[x,y]$ is central, then $e^xe^y=e^{x+y+c}=e^{c}e^{x+y}$ by the BCH formula.
In particular, for
\eq{b=\sum_{j=1}^s t_jb_j,\qquad
c=\oh\sum_{i<j}t_it_j[b_i,b_j],}
we obtain
\eq{w:=e^{t_1b_1}\dots e^{t_s b_s}=e^c e^b.}
Therefore
\eq{\ta_w(f):=\ta(fw)/\ta(w)=\ta(fe^b)/\ta(e^b).}

\begin{theorem}[Twisted Wick formula]
\label{wick twisted:th}
For $a_1,\dots,a_r\in\Ho$, we have $\ta_w(a_i)=\ta(a_ib)$ and
\eq{\label{wick twisted}
\ta_w\rbr{a_1\dots a_r}
=\sum_{\pi\in\M_r}
{\prod_{\set{i<j}\in\pi}\ta(a_ia_j)}
\prod_{\set i\in\pi}\ta(a_ib).
}
\end{theorem}

\begin{proof}
By \cref{wick-exp1}, we have
\begin{multline*}
\ta\rbr{a_1\dots a_r w}
=e^c\ta\rbr{a_1\dots a_r e^b}
=e^c\ta(e^b)\sum_{\pi\in\M_r}
\prod_{\set{i<j}\in\pi}\ta(a_i a_j)
\prod_{\set i\in\pi} \ta(a_i b)\\
=\ta(w)\sum_{\pi\in\M_r}
\prod_{\set{i<j}\in\pi}\ta(a_i a_j)
\prod_{\set i\in\pi} \ta(a_i b).
\end{multline*}
In particular, we obtain $\ta_w(a_i)=\ta(a_ib)$.
\end{proof}

\begin{lemma}
For $b=\sum_{j=1}^s t_jb_j$ and $w=\prod_{j=1}^s e^{t_jb_j}$, we have
\begin{equation}\label{ta(w)}
\ta(e^b)=\exp\rbr{\oh\ta(b^2)},\qquad
\ta\rbr{w}
=\exp\rbr{
\frac12\sum_{j=1}^s t_j^2\ta(b_j^2)
+\sum_{i<j}t_it_j\ta(b_ib_j)
}.
\end{equation}
\end{lemma}
\begin{proof}
By the ordinary Wick formula, we have $\ta(e^b)=\sum_{k\ge0}\frac{1}{(2k)!}\frac{(2k)!}{2^kk!}\ta(b^2)^{k}
=\exp\rbr{\oh\ta(b^2)}$.
For $c=\oh\sum_{i<j}t_it_j[b_i,b_j]$, we have
\[c+\oh\ta(b^2)
=\oh\sum_jt_j^2\ta(b_j^2)
+\oh\sum_{i<j}t_it_j\ta([b_i,b_j]+b_ib_j+b_jb_i)
.\]
Therefore
$\ta(w)=e^c\ta(e^b)=\exp\rbr{c+\oh\ta(b^2)}
=\exp\rbr{\oh\sum_jt_j^2\ta(b_j^2)
+\sum_{i<j}t_it_j\ta(b_ib_j)}$.
\end{proof}

For $u,v\in H$, let $\phi_u(z)=\sum_{n\ne0}u_{-n}\frac{z^n}n$ and let
\begin{gather}
w=W_{u,v}(z)=\nop{e^{\vi_{u,v}(z)}}=e^{\vi_u(z)_+}e^{\vi_v(z)_-},
\label{W_uv}
\\
\vi_{u,v}(z)=\vi_u(z)_++\vi_v(z)_-
=\sum_{n>0}\frac{u_{-n}z^n-v_nz^{-n}}n.
\end{gather}
For $a\in H$, we have
\eq{\label{ta_w(a_n)}
\begin{gathered}
\ta_w(\fa_n(a))
=\ta(\fa_n(a)\fa_{-n}(u))\frac{z^n}n
=(a,u)\frac{z^n}{1-q^n},\qquad n>0,
\\
\ta_w(\fa_{n}(a))
=\ta(\fa_{n}(a)\fa_{-n}(v))\frac{z^{n}}{n}
=(a,v)\frac{z^n}{1-q^n},\qquad n<0.
\end{gathered}
}

\begin{lemma}\label{ta(W)}
For $W_{u,v}(z)=e^{\vi_u(z)_+}e^{\vi_v(z)_-}$, we have
\[\ta(W_{u,v})=\ph q^{(u,v)}.\]
In particular, $\ta(W_{u,v})=1$ if $(u,v)=0$.
\end{lemma}
\begin{proof}
Applying \cref{ta(w)} to $w=e^{\vi_u(z)_+}e^{\vi_v(z)_-}$
and using $\ta(\vi_u(z)_+^2)=\ta(\vi_v(z)_-^2)=0$, we obtain
\begin{multline}
\ta(W_{u,v})
=\exp\rbr{\ta(\vi_u(z)_+\vi_v(z)_-)}
=\exp\rbigg{-\sum_{n>0}\frac{\ta(u_{-n}v_n)}{n^2}}
\\
=\exp\rbigg{-(u,v)\sum_{n>0}\frac1n\frac{q^n}{1-q^{n}}}
=\ph q^{(u,v)}.
\tag*{\qedhere}
\end{multline}
\end{proof}

\section{Okounkov's conjecture}
\label{ok:sec}
We now combine the preceding results to establish
that normalized $q$-brackets are multiple $q$-zeta values.
The proof reduces the geometric generating series to graph sums whose structure is constrained by the cohomology of the surface.
Their arithmetic properties can then be studied using the Fourier and summation techniques developed above.

\subsection{Formulation of the theorem}
Let $X$ be a smooth projective surface.
For simplicity we assume that $H^\odd(X,\bQ)=0$.
Let $E_1,\dots,E_r$ be vector bundles over $X$
and let $k_1,\dots,k_r\ge0$.
Consider the $q$-bracket
\[
\ang{\ch_{k_1}^{E_1}\dots \ch_{k_r}^{E_r}}
=\sum_{n\ge0}q^n\int_{\hb Xn}
\ch_{k_1}(\hb{E_1}n)\dots \ch_{k_r}(\hb{E_r}n)
\cdot \fe_t(T_{\hb Xn})
\]
and the normalized $q$-bracket $\ang{f}'=\ang{f}/\ang{}$.
We will prove the following result.

\begin{theorem}\label{Ok1}
We have
\[
\ang{\ch_{k_1}^{E_1}\dots \ch_{k_r}^{E_r}}'\in\Zq[t].
\]
\end{theorem}

First, let us interpret the above $q$-bracket as a $q$-trace
using the results of \cite{carlsson_vertex,carlsson_exts}
(see also \cite{hongminh_quasia} for more details).
Given a line bundle $L$ with $a=c_1(L)$, let
\[W(L,z)
=W_{a-K,a}(z)
=\exp\rbr{\phi_{a-K}(z)_+}\exp\rbr{\phi_{a}(z)_-},
\]
where $\phi_a(z)=\sum_{n\ne 0}a_{-n}\frac{z^n}n$.

\begin{theorem}[Carlsson-Okounkov]
\label{CO:th}
Let $\Phi=(\Phi(n))_{n\ge0}\in\prod_{n\ge0}\bH_n$
and let $\Phi_\bul\in\End(\bH)$ be the corresponding cup product operator.
Then
\[\sum_{n\ge0}q^n\int_{\hb Xn}\Phi(n)\cdot\fe_t(T_{\hb Xn})
=\tr_q\rbr{\Phi_\bul W(L,z)}
\]
for the equivariant line bundle $L=\cO_X\ts e^t$.
\end{theorem}

For $L=\cO_X\ts e^t$, we have $c_1(L)=t$
and
\[w:=W(L,z)=W_{t-K,t}(z).\]
Since $(t-K,t)=0$, we have $\ta(w)=1$ by \cref{ta(W)}.

For $\Phi(n)=\ch_{k_1}(\hb{E_1}n)\dots \ch_{k_r}(\hb{E_r}n)$,
we have $\Phi_\bul=\fch_{k_1}(E_1)\dots\fch_{k_r}(E_r)$
and
\eq{
\ang{\ch_{k_1}^{E_1}\dots \ch_{k_r}^{E_r}}'
=\frac{\tr_q(\Phi_\bul w)}{\tr_q(w)}
=\frac{\ta(\Phi_\bul w)}{\ta(w)}=\ta_w(\Phi_\bul),
}
where $\ta(f)=\tr_q(f)/\tr_q(1)$ and
$\ta_w(f)=\ta(fw)/\ta(w)=\tr_q(fw)/\tr_q(w)$.

Instead of the Chern character operators $\fch_{k_i}(E_i)$,
we may consider the tautological operators
$\fG_{k_i}(a_i)$ for $a_i\in H=H^*(X,\bQ)$.
Then \cref{Ok1} admits the following generalization.

\begin{theorem}\label{Ok2}
Let $w=W_{t-K,t}(z)$, $a_1,\dots,a_N\in H$ and
$j_1,\dots,j_N\ge0$.
Then
\eq{\label{prod ta_w}
\ta_w\rbr{\prods_{k=1}^N\fG_{j_k}(a_k)}\in\Zq[t].}
\end{theorem}

\subsection{First reduction}
Let $w=W_{u,v}(z)$ \cref{W_uv} with $u=t-K$ and $v=t$.
Let $\cE=\set{1,e,K,K^2}$.
For $\ep\in\cE$, let $\kp(\ep)\in\set{0,1,2}$ be the
power of $K$ in $\ep$, and let $\n\ep$ denote its
cohomological degree.
Then
\[
\fG_j(a)
=(-1)^j
\sum_{\ep\in\cE}
\sum_{\ov{\n\la=0}{\ell(\la)=j+2-\n\ep/2}}
\frac{g_\ep(\la)}{\la!}\fa_\la(\ep a),
\]
where $g_\ep(\la)$ were determined in \cref{cor:g_ep}.
This result implies that \cref{prod ta_w}
is a linear combination of expressions
$\ta_w(\cT_1\cdots\cT_N)$, where
\eq{
\cT_k
=\sum_{\ov{\bn_k\in\Zs^{r_k}}{\n{\bn_k}=0}}
P_k(\bn_k)\rbr{\prods_{I\in\cI_k}\n{n_I}}
\fa_{\bn_k}(b_k),
}
with
\[
\ep_k\in\cE,\qquad
r_k=j_k+2-\n{\ep_k}/2,\qquad
P_k\in\bQ[n_1,\dots,n_{r_k}],\qquad b_k=\ep_k a_k,
\]
$\cI_k\sbs 2^{[r_k]}$ consisting of $\le\kp(\ep_k)$ disjoint subsets,
and $n_I=\sum_{i\in I}n_{k,i}$ for $I\sbs[r_k]$.

Let
\[
V=[N],\qquad
F=\bigsqcup_{k\in V}\set{k}\xx[r_k],
\]
and let $p:F\to V$ be the projection
and $F_k=p\i(k)\iso [r_k]$ for $k\in V$.
We consider the lexicographic total order on $F$.
Let
\[\cN=\sets{\bn\in\Zs^F}{\n{\bn_k}
=\sums_{i\in F_k}n_i=0\ \forall k\in V},\]
where $\bn_k=\bn|_{F_k}$.
Then
\begin{gather}
\ta_w(\cT_1\cdots\cT_N)=\sum_{\bn\in\cN}
P_V(\bn) \ta_w\rbr{\prods_{k=1}^N\fa_{\bn_k}(b_k)},
\label{ta_w:cT}
\\
P_C(\bn)
=\prods_{k\in C}\rbr{P_k(\bn_k)\prods_{I\in\cI_k}\n{n_I}},
\qquad C\sbs V.
\end{gather}

We will apply twisted Wick \cref{wick twisted:th}
to express $\ta_w\rbr{\prod_{k=1}^N\fa_{\bn_k}(b_k)}$
as a sum over matchings $\pi\in\cM_F$.
For $\pi\in\cM_F$, let $\Ga_\pi=(V,F,p,\pi)$
be the graph with legs (see \cref{graph w legs})
and let
\eq{\label{N_pi}
\cN_\pi=\cN_{\Ga_\pi}
=\sets{\bn\in\Zs^F}{
\n{\bn_k}=0\ \forall k\in V,\
n_i+n_j=0\ \forall\set{i,j}\in E(\Ga_\pi)}.
}
For a subgraph $C\sbs \Ga_\pi$ and $\bn\in\cN$,
let
\[r_C(\bn)=\#\sets{i\in L(C)}{n_i>0},\qquad
b_C=\prod_{k\in V(C)}b_k\in H.
\]

\begin{theorem}
For $\bn\in\cN$, we have
\begin{equation}\label{ta-A}
\ta_w\rbr{\prods_{k=1}^N\fa_{\bn_k}(b_k)}
=\sum_{\pi\in\M_F}A_\pi(\bn),
\end{equation}
where $A_\pi(\bn)=0$ unless $\bn\in\cN_\pi$, in which case
\[A_\pi(\bn)
=\prod_{C\in\pi_0(\Ga_\pi)}f_C(\bn)\eps_C(\bn),\qquad
\eps_C(\bn)=\eps\rbr{e^{b_1(C)}b_C u^{r_C(\bn)}v^{\n{L(C)}-r_C(\bn)}},
\]
\[
f_{C}(\bn)
=\prod_{\ov{\set{i<j}\in E(C)}{p(i)<p(j)}}\rb(n_i)
\prod_{\ov{\set{i<j}\in E(C)}{p(i)=p(j)}}\rb(-\n{n_i})
\prod_{i\in L(C)}\ra(n_i).
\]

\end{theorem}
\begin{proof}
Let $\de_{r_k}(b_k)=\ts_{i\in F_k}x_i$ (we suppress the summation sign).
Then
\[
\prod_{k=1}^N\fa_{\bn_k}(b_k)
=\prod_{k=1}^N \nop{\prod_{i\in F_k}\fa_{n_i}\rbr{x_i}}.
\]

Applying $\ta_w$ and using the twisted Wick formula
\cref{wick twisted} for $w=W_{u,v}(z)=e^{\phi_u(z)_+}e^{\phi_v(z)_-}$
and $b=\phi_u(z)_++\phi_v(z)_-$,
we obtain a sum over $\pi\in\cM_F$.
For a pair $\set{i<j}\in \pip$,
the contraction vanishes unless $n_i+n_j=0$.
Under this condition it is \cref{ta ab}
\[
\begin{cases}
\ta\rbr{\fa_{n_i}(x_i)\fa_{n_j}(x_j)}
=(x_i,x_j)\rb(n_i),
& p(i)<p(j),\\
\ta\rbr{\nop{\fa_{n_i}(x_i)\fa_{n_j}(x_j)}}
=(x_i,x_j)\rb(-\n{n_i}), & p(i)=p(j).
\end{cases}
\]
Indeed, operators in different blocks $F_k$ retain their block
order, whereas operators within a block are normally ordered.
For $\set i\in \pis$, the contribution is \cref{ta_w(a_n)}
\[
\ta_w\rbr{\fa_{n_i}(x_i)}
=\frac{z^{n_i}}{1-q^{n_i}}
\begin{cases}
(x_i,u),&n_i>0,\\
(x_i,v),&n_i<0.
\end{cases}
\]
Since $\bn\in\cN_\pi$, we have
$\sum_{i\in L(\Ga_\pi)}n_i=\sum_{i\in F}n_i-\sum_{\set{i,j}\in E(\Ga_\pi)}(n_i+n_j)=0$.
Hence the factors $z^{n_i}$ cancel.
The scalar factors from the pairs and singletons therefore
give precisely $\prod_{C\in\pi_0(\Ga_\pi)}f_{C}(\bn)$.

It remains to take the product of the Frobenius contractions.
For each connected component $C\in\pi_0(\Ga_\pi)$, the vertices carry
the classes $b_k$, the edges contract tensor factors using
the Frobenius pairing, and the unmatched flags are paired
with $u$ or $v$ according to the sign of their modes.
By \cref{lem:frobenius-graph}, this contraction equals
\[
\eps\rbr{e^{b_1(C)}b_C u^{r_C(\bn)}v^{\n{L(C)}-r_C(\bn)}}.
\]
Multiplying over the connected components proves the formula.
\end{proof}

Since $P_V(\bn)=\prod_{\Ga\in\pi_0(\Ga_\pi)}P_\Ga(\bn)$ for $\pi\in\M_F$, we obtain
\eq{
\ta_w(\cT_1\dots\cT_N)
=\sum_{\pi\in\M_F}\sum_{\bn\in\cN_\pi}P_V(\bn)A_\pi(\bn)
=\sum_{\pi\in\M_F}\prod_{\Ga\in\pi_0(\Ga_\pi)}
\rbr{\sum_{\bn\in\cN_\Ga}P_\Ga(\bn)f_\Ga(\bn)\eps_\Ga(\bn)},
}
where $\cN_\Ga\sbs\Zs^{F(\Ga)}$ was defined in \cref{N_pi}.
This implies that in order to prove \cref{Ok2}
it suffices to show
that
\eq{\label{Ga contrib}
\sum_{\bn\in\cN_\Ga}P_\Ga(\bn)f_\Ga(\bn)\eps_\Ga(\bn)\in \Zq[t]
}
for every
connected component $\Ga\in\pi_0(\Ga_\pi)$ and $\pi\in\M_F$.
Therefore we can assume that $\pi\in\M_F$
is a matching such that the graph with legs
\[\Ga=(V,F,p,\pi)\]
is connected.
Let $E=E(\Ga)$ and let $L=L(\Ga)$.

\begin{lemma}
Assume that $\Ga$ is connected and $b_1(\Ga)\ge1$.
Then
\[\sum_{\bn\in\cN_\Ga}P_\Ga(\bn)f_\Ga(\bn)\eps_\Ga(\bn)\in\QM[t].\]
\end{lemma}
\begin{proof}
Assume that
$\eps_\Ga(\bn)
=\eps\rbr{e^{b_1(\Ga)}b_\Ga u^{r_\Ga(\bn)}v^{\n{L}-r_\Ga(\bn)}}\ne0$ for some $\bn$.
Since $e^2=0$, we obtain $b_1(\Ga)=1$.
Since $u=t-K$ and $v=t$, we obtain $\eps_\Ga(\bn)=t^{\n L}\eps(e b_\Ga)$
and we need to show that $\vi_\Ga=\sum_{\bn\in\cN_\Ga}P_\Ga(\bn)f_\Ga(\bn)\in\QM$.
Degree considerations imply that $b_\Ga$ has degree zero, hence $\ep_k=1$ for all $k\in V$.
Therefore $\cI_k=\es$ and
\[
P_\Ga(\bn)=\prods_{k\in V}P_k(\bn_k)\in
\bQ[n_i\col i\in F].
\]

By linearity, we may assume that $P_\Ga$ is a monomial.
Suppose first that $\Ga$ contains a loop $\set{i<j}$.
Since $n_j=-n_i$, its mode does not contribute to the vertex conservation equations.
The summation over this mode therefore separates from the remaining modes.
We have
\[
\sum_{n\ne0}n^k\rb(-\n n)
=\sum_{n>0}(1+(-1)^k)n^{k+1}\frac{q^n}{1-q^n}
=\begin{cases}
0&\text{odd }k,\\
2\sum_{n>0}Q_{k+2}(q^n)&\text{even }k.
\end{cases}
\]
This expression belongs to $\QM$. We may therefore remove the loop.

We associate with $\Ga$
the acyclic quiver $Q$ with the set of vertices $V\sqcup\set*$
and the arrows $a_{ij}:p(j)\to p(i)$ for $\set{i<j}\in E$ and
$a_i:*\to p(i)$ for $i\in L$.
Let $\dd:\bZ^{Q_1}\to\bZ^{Q_0}$ be the divergence \cref{boundary}.
There is a bijection between the set $\cN_\Ga$
and the set $\sets{\bm\in\Zs^{Q_1}}{\dd(\bm)=0}$,
where $\bn\in\cN_\Ga$ corresponds to $\bm\in\Zs^{Q_1}$
with $m_{a_{ij}}=n_i=-n_j$ for $\set{i<j}\in E$
and $m_{a_i}=n_i$ for $i\in L$.

For $\set{i<j}\in E$ and $a=a_{ij}$, we have
$\rb(n_i)=\what f_a(m_a)$, where $f_a(z)=\Phi_2(z)\in\cC_1$.
For $i\in L$ and $a=a_i$,
we have $\ra(n_i)=\what f_a(m_a)$,
where $f_a(z)=\Phi_1(z)\in\cC_1$.
Then
\[\vi_\Ga=\sum_{\ov{\bm\in\Zs^{Q_1}}{\dd \bm=0}}
P_\Ga(\bm)\prod_{a\in Q_1}\what f_a(\bm_a)\in\QM
\]
by \cref{QM acyclic}.
\end{proof}

\subsection{Reduction to tree sums}
It remains to assume that the graph with legs $\Ga=(V,F,p,\pi)$
is connected and $b_1(\Ga)=0$, hence $\Ga$ is a tree.
Then the contribution of $\Ga$ is \cref{Ga contrib}
\[
A(q)=\sum_{\bn\in\cN_\Ga}
P_\Ga(\bn)f_\Ga(\bn)
\eps\rbr{b_\Ga(t-K)^{r(\bn)}t^{\n L-r(\bn)}},
\]
where $r(\bn)=\#\sets{i\in L}{n_i>0}$
and
\[
P_\Ga(\bn)=P(\bn)\prod_{k\in V}\prod_{I\in\cI_k}\n{n_I},\qquad
f_{\Ga}(\bn)
=\prod_{\set{i<j}\in E}\rb(n_i)\prod_{i\in L}\rho(n_i),
\]
for a polynomial $P(\bn)$.
We have
\[
\eps\rbr{b_\Ga(t-K)^{r(\bn)}t^{|L|-r(\bn)}}
=
\sum_{j=0}^{\min(2,|L|)}
(-1)^jt^{|L|-j}
\binom{r(\bn)}j\eps(b_\Ga K^j).
\]
For $S\sbs L$, define
\[
A_{S}(q)
=\sum_{\bn\in\cN_\Ga}
P_\Ga(\bn)f_\Ga(\bn)\prod_{i\in S}H(n_i),
\]
where $H(x)=\one_{x\ge0}$ is the Heaviside step function.
Then
\[
A(q)
=\sum_{j=0}^{\min(2,|L|)}
(-1)^jt^{|L|-j}\eps(b_\Ga K^j)
\sum_{\ov{S\sbs L}{\n S=j}}
A_{S}(q).
\]

Let $\cI=\bigsqcup_{k\in V}\cI_k$.
Recall that $\n{\cI_k}\le\kp(\ep_k)$, where $\kp(\ep_k)$ is the power of $K$ in $\ep_k\in\set{1,e,K,K^2}$.
If $\eps(b_\Ga K^j)\ne0$, then
\[
2\n\cI+2j
\le 2\sum_{k\in V}\kp(\ep_k)+2j
\le \deg b_\Ga+2j=4,
\]
hence $\n \cI+j\le2$.
Therefore it suffices to prove that
$A_S(q)\in\Zq$
whenever $\n\cI+\n S\le 2$.

Since $\rb(n_i)=n_i\ra(n_i)$,
we can absorb $n_i$ into $P$ and assume that
\[
f_{\Ga}(\bn)
=\prod_{\set{i<j}\in E}\rho(n_i)\prod_{i\in L}\rho(n_i).
\]

We associate with $\Ga$ the quiver $Q$ with the set of vertices $V\sqcup\set*$,
the internal arrows $a_{ij}:p(j)\to p(i)$ for $\set{i<j}\in E$ and
the framing arrows $a_i:*\to p(i)$ for $i\in L$.
An element $\bn\in\cN_\Ga$ corresponds to $\bm\in\Zs^{Q_1}$
satisfying $\dd(\bm)=0$ \cref{boundary},
where $m_{a_{ij}}=n_i=-n_j$ for $\set{i<j}\in E$
and $m_{a_i}=n_i$ for $i\in L$.

We identify $F_k$ with the set of arrows incident with $k$.
For $I\sbs F_k$, define
\[\dd_I(\bm)=\sum_{\ov{a\in I}{t(a)=k}}m_a
-\sum_{\ov{a\in I}{s(a)=k}}m_a.
\]
Then $n_I=\dd_I(\bm)$ for $\bm\in\Zs^{Q_1}$ corresponding to
$\bn\in\cN_\Ga$.
Identifying $S\sbs L$ with the corresponding set of framing arrows, we obtain
\[
A_{S}(q)
=\sum_{\ov{\bm\in\Zs^{Q_1}}{\dd(\bm)=0}}
P(\bm)\prod_{I\in\cI}|\dd_I(\bm)|
\prod_{a\in Q_1}\rho(m_a)\prod_{a\in S}H(m_a),
\]
where $P(\bm)$ is a polynomial.

We now simplify the factors $\n{\dd_I(\bm)}$ by splitting vertices.
If $\cI_k=\set{I,J}$ and $I\sqcup J=F_k$,
the conservation gives
$\n{\dd_I(\bm)}\n{\dd_J(\bm)}=\dd_I(\bm)^2$
and we can absorb this factor into $P$.
If $F_k\in\cI_k$, then $A_S(q)=0$.
Therefore we can assume that $R_k=F_k\ms \bigcup_{I\in\cI_k}\ne\es$
for all $k$.

We replace every vertex $k\in V$ by $(k,I)$
for $I\in\cI'_k=\cI_k\sqcup\set {R_k}$.
Every original arrow $a:k\to l$
is replaced by
$a:(k,I)\to (l,J)$, where $I\in\cI'_k$ and $J\in\cI'_l$
are the unique parts containing $a$.
Every arrow $a:*\to k$ is replaced by $a:*\to (k,I)$, where $a\in I$.
For every $I\in\cI_k$
consider an auxiliary arrow $a_I:(k,I)\to (k,R_k)$.

Let $Q'$ be the new quiver.
Let $\bm\in\Zs^{Q_1}$ and let $\ell_I\in\Zs$ correspond to the auxiliary arrow $a_I$ for $I\in\cI$.
Conservation at the new vertices gives
\[\dd_I(\bm)-\ell_I=0,\qquad \dd_{R_k}(\bm)+\sum_{I\in\cI_k}\ell_I=0.\]
These equations are equivalent to
the original conservation $\dd(\bm)=0$, together with $\ell_I=\dd_I(\bm)$ for every $I\in\cI$.
The internal graph of the new quiver (obtained by removing $*$) is again a tree.
Let $\cA\sbs Q'_1$ be the set of the original arrows.
Then $A_S(q)$ can be written as
\[\sum_{\ov{\bm\in\Zs^{Q'_1}}{\dd(\bm)=0}}
P(\bm)\prod_{I\in\cI}\n{\ell_I}
\prod_{a\in\cA}\rho(m_a)\prod_{a\in S}H(m_a),
\]
where $\ell_I=m_{a_I}$.
We have $\n{x}=x\sgn(x)$ and $\sgn(x)=2H(x)-1$ for $x\ne0$.
Therefore we can absorb $\prod_{I\in\cI}\ell_I$ into $P$
and expand the remaining sign factors.
Therefore it suffices to prove that
\[\sum_{\ov{\bm\in\Zs^{Q'_1}}{\dd(\bm)=0}}
P(\bm)\prod_{a\in\cA}\rho(m_a)\prod_{a\in S'}H(m_a)
\in\Zq
\]
for every $S'\sbs Q'_1$
satisfying $\n{S'}\le\n\cI+\n S\le 2$.

Finally, since $1=\rho(n)+\rho(-n)$,
we can introduce factors $\rho(\pm m_a)$ for all auxiliary arrows $a\in Q'_1\ms \cA$.
Thus it suffices to prove that
\eq{\label{BD red-5}
\sum_{\ov{\bm\in\Zs^{Q'_1}}{\dd(\bm)=0}}
P(\bm)\prod_{a\in Q'_1}\rho(\eps_a m_a)\prod_{a\in S'}H(m_a)
\in\Zq,
}
where $\n {S'}\le 2$, $\eps_a\in\set{\pm1}$ for $a\in Q'_1$,
and $\eps_a=1$ for all $a\in Q'_1$ incident with $*$.

\begin{remark}
We prove the above statement in the next section.
A simpler proof would follow from \cref{conj1}.
Indeed, in \cref{BD red-5} we can choose $\bm\in\Zs^{Q'_1}$ with fixed signs $\si_a=\sgn(m_a)$ and
define the quiver $Q$ by reversing the arrows $a$ with $\si_a=-1$.
Then $n_a=\si_a m_a>0$ and $\eps_am_a=\eps_a\si_a n_a$.
Let $Q''$ be the quiver obtained by reversing the arrows with $\eps_a\si_a=-1$ in $Q$ or, equivalently,
reversing the arrows with $\eps_a=-1$ in $Q'$.
The arrows with $\eps_a=-1$ are contained in the internal tree,
hence $Q''$ is acyclic.
Therefore the sum in \cref{BD red-5} corresponding to the fixed signs $\si_a$ is
\[
\sum_{\ov{\bn\in\Ns^{Q_1}}{\dd_Q(\bn)=0}}P'(\bn)
\prod_{a\in Q_1}\rho(\eps_a\si_a n_a)
\]
which belongs to $\BD$ by \cref{conj1}.
\end{remark}

\subsection{Framed tree theorem}
We conclude by proving a summation theorem for framed trees, which provides the final ingredient in the proof of our version of Okounkov's conjecture.
The restriction to at most two marked arrows makes it possible to reduce these graph sums to the bi-bracket summation results established earlier.

Let $Q'=(V',E',s,t)$ be a quiver with a fixed vertex $*\in V'$
such that for $V=V'\ms\set*$
the internal graph $Q'|_V=(V,E)$ is a tree and all arrows in $E_*=E'\ms E$ are of the form $*\to i$ for $i\in V$ (we call them framing arrows).
We define
\[
\dd:\bZ^{E'}\to\bZ^V,\qquad
\dd_i(\bn)=\sum_{t(e)=i}n_e-\sum_{s(e)=i}n_e.
\]
For $i\in V$ and $S\sbs E'$, let
\[
\dd_i^S(\bn)=
\sum_{\substack{e\in S\\t(e)=i}}n_e
-\sum_{\substack{e\in S\\s(e)=i}}n_e.
\]

Recall that $\rho(n)=\hPhi_1(n)=\frac1{1-q^n}$ for $n\ne0$
and $\rho(n)+\rho(-n)=1$.

\begin{theorem}
\label{thm:framed-tree}
Let $\eps\in\set{\pm1}^{E'}$ with $\eps_e=1$ for all $e\in E_*$,
and let $\Sa\sbs E'$ with $\n\Sa\le2$.
Then, for every polynomial
$P\in\bQ[n_e\col e\in E']$,
\eq{\label{BD1}
\sum_{\substack{\bn\in\Zs^{E'}\\ \dd(\bn)=0}}
P(\bn)\prod_{e\in E'}\rho(\eps_en_e)
\prod_{e\in \Sa}H(n_e)
\in\Zq.
}
\end{theorem}

\begin{proof}
\clm{Reduction to $\eps=1$}
We reverse every arrow $e$ with $\eps_e=-1$ and replace
$n_e$ by $m_e=\eps_en_e$ for all $e\in E'$.
This preserves the conservation equations.
All framing arrows remain unchanged, and the internal
graph is still a tree.
Moreover,
\[
\rho(\eps_en_e)=\rho(m_e),\qquad
H(n_e)=
\begin{cases}
H(m_e),&\eps_e=1,\\
1-H(m_e),&\eps_e=-1.
\end{cases}
\]
Expanding the Heaviside factors reduces the theorem
to the case $\eps_e=1$ for every $e\in E'$.
We order the vertices so that $i<j$ for every arrow $i\to j$.

\clm{Properness}
Let $\pp n=\max\set{0,n}$ and $\pn n=\max\set{0,-n}$,
so that
\[n=\pp n-\pn n,\qquad \n n=\pp n+\pn n.\]
We have $\rho(n)=\sum_{k\ge0}q^{kn}$ and
$\rho(-n)=1-\rho(n)=-\sum_{k\ge1}q^{kn}$ for $n>0$.
Therefore $\ord_q\rho(n)=\pn n$
and the valuation of $\prod_{e\in E'}\rho(n_e)$ is
\[
v(\bn)=\sum_{e\in E'}\pn{n_e}.
\]
Since $\sum_{e\in E_*}n_e=0$,
we have
\[
v(\bn)\ge\sum_{e\in E_*}\pn{n_e}
=\sum_{e\in E_*}\pp{n_e},
\]
hence $\sum_{e\in E_*}\n{n_e}\le2v(\bn)$.
For an internal arrow $e\in E$,
the conservation property $\dd(\bn)=0$
implies that
$\n{n_e}\le\sum_{f\in E_*}\n{n_f}\le 2v(\bn)$.
This implies that there are only finitely many summands in \cref{BD1} with a fixed valuation $v(\bn)$.
Therefore the whole sum is $q$-proper.

\clm{Reduction}
We can assume that $P(\bn)=\prod_{e\in E'}n_e^{p_e}$ for some $p_e\ge0$.
Because of the factors $H(n_e)$ in~\cref{BD1},
we can assume that $d_e=n_e>0$ for $e\in\Sa$.
Then \cref{BD1} has the form
\eq{\label{BD3}
A=\sum_{\bd\in\pN^\Sa}
\rbr{\prod_{e\in \Sa}d_e^{p_e}\rho(d_e)}\Phi(\bd),\qquad
\Phi(\bd)=
\sum_{\ov{\bn\in\Zs^{E'\ms \Sa}}{\dd(\bn,\bd)=0}}
\prod_{e\in E'\ms \Sa}n_e^{p_e}\rho( n_e).
}
Deleting the arrows from $\Sa\cap E$
splits the internal tree $(V,E)$ into components
$C_1,\ldots,C_s$.
Let $M\sbs V$ be the set of vertices incident with $\Sa$,
and let $M_\al=M\cap C_\al$.
An internal arrow in $\Sa$ has its endpoints in distinct components, while a framing arrow has only one endpoint in $V$.
Therefore $\n{M_\al}\le\n{\Sa}\le2$.
The conservation equation $\dd(\bn,\bd)=0$ is equivalent to
\[
\dd_i^{E'\setminus \Sa}(\bn)=-\dd_i^\Sa(\bd)
=:m_i(\bd)
\qquad \forall i\in V.
\]
We have $m_i(\bd)=0$ for $i\in V\ms M$.
Let $m^\al(\bd)=(m_i(\bd))_{i\in C_\al}$ for $\al\in[s]$.

Consider the variables $z_i$ for $i\in V$, and let $z_*=1$.
For every arrow $e:i\to j$, let $z_e={z_j}/{z_i}$
and let
\[
U_e(z):=\Phi_{p_e+1}(z)
=D_{z}^{p_e}\rbr{\Phi_1(z)}
=\sum_{n\ne0}n^{p_e}\rho(n)z^{n}\in\cC_1.
\]
For every component $C_\al$, define
\[G_\al=\prod_{\ov{e\in E'\ms \Sa}{t(e)\in C_\al}}U_e(z_e)
\in\cC_{\n{C_\al}}.
\]
The variables in the products $G_\al$ are disjoint.
Moreover, the exponent of $z_i$ is precisely
$\dd_i^{E'\ms\Sa}(\bn)$.
Hence coefficient extraction gives
\[\Phi(\bd)
=\prod_{\al=1}^s \what G_\al(m^\al(\bd)).
\]

\clm{Fourier--Lambert summation}
Let us prove that \cref{BD3} belongs to $\Zq$.
If $\Sa=\es$, this follows directly from the Constant term theorem.

Let $\Sa=\set{e_1,\dots,e_r}$ for $1\le r\le2$, and let $d_k=d_{e_k}$.
By \cref{fourier}, on every chamber determined by the signs of the subset sums of $m^\al(\bd)$
(meaning $L(\bd)=\sum_{i\in I}m_i(\bd)$ for some $I\sbs C_\al$),
every $\what G_\al(m^\al(\bd))$
is a finite $\QM$-linear combination
of polynomials in $\bd$ times factors $f(q^{\n{L(\bd)}})$,
where $f\in\RR$ and $L$ is a nonzero subset sum of $m^\al(\bd)$.

Since $\n {M_\al}\le 2$, the possible forms $L$ are $m_i$ and $m_i+m_j$ for $i,j\in M_\al$.
For $r=1$ they are, up to sign, equal to $d_1$.
For $r=2$ they are, up to sign, equal to
\[\set{d_1,d_2,d_1+d_2,d_1-d_2}.\]
The internal tree has at most one component incident with both arrows from $\Sa$.
Therefore at most one of $d_1+d_2$ and $d_1-d_2$ occurs, up to sign.
If $d_1-d_2$ does not occur, we take $(u,v)=(d_1,d_2)$.
If $d_1-d_2$ occurs, we split the summation domain $\pN^2$ in \cref{BD3} into regions $d_1>d_2$, $d_1<d_2$ and $d_1=d_2$, using coordinates
\[(u,v)=(d_1-d_2,d_2),\qquad (u,v)=(d_2-d_1,d_1)\]
on the first two regions.
On every two-dimensional region,
the expression $\Phi(\bd)$
is a finite $\QM$-linear combination of
\[P(u,v)f_1(q^u)f_2(q^v)f_3(q^{u+v}),\]
where $P$ is a polynomial and $f_i\in\RR$.
The remaining factor
$\prod_{e\in \Sa}d_e^{p_e}\rho(d_e)$
is of the same form, since $\rho(d)=1+Q_1(q^d)$,
where $Q_1(q)=\frac q{1-q}\in\RR$.
On the one-dimensional regions,
the expression $\Phi(\bd)$ has the analogous form
$P(u)f(q^u)$, where $P$ is a polynomial and $f\in\RR$.
The remaining factor $\prod_{e\in \Sa}d_e^{p_e}\rho(d_e)$
has the same form.

The summand in \cref{BD3} is $q$-proper by the first part of the proof.
This property is preserved under restriction to the regions above
and under the corresponding changes of variables.
On each region, the summand satisfies the assumptions of \cref{proper sum}, with the nested family $\set{\set1,\set2,\set{1,2}}$ in dimension two and $\set{\set1}$ in dimension one.
Hence the sum over each region belongs to $\Zq$,
and therefore $A\in\Zq$.
\end{proof}

\bibliographystyle{halpha}
\bibliography{qmzv}
\end{document}